\documentclass[1p, 11p]{elsarticle}

\usepackage{graphicx}
\usepackage{caption}
\usepackage{subcaption}
\usepackage{comment}
\usepackage{amssymb}
\usepackage{amsmath}
\usepackage{dirtytalk}
\usepackage{algorithm}
\usepackage{algorithmic}
\biboptions{numbers,sort&compress}
\usepackage{lineno}

\usepackage{xcolor}
\usepackage{hyperref}
\usepackage{siunitx}
\usepackage{tikz}
\usepackage{pgfplots, pgfplotstable}
\pgfplotsset{compat = newest}
\usetikzlibrary{matrix}
\usetikzlibrary{plotmarks}
\pgfplotsset{compat=newest}
\usepackage{makecell}
\usepackage[export]{adjustbox}
\usepackage{color, colortbl}

\definecolor{Gray}{gray}{0.8}

\newlength{\offsetpage}
{\end{center}%
\end{adjustwidth}}
\newcommand\commentout[1]{\relax}

\journal{Engineering Fracture Mechanics}

\begin{document}

\begin{frontmatter}

\title{Quasistatic cohesive fracture with an alternating direction method of multipliers}

\author[appliedmath]{James I. Petrie}
\author[mechanicalUW,hmi]{M. Reza Hirmand}
\author[yorktesserea]{Katerina D. Papoulia}
\cortext[cor]{Corresponding author:}
\ead{papouli@yorku.ca}
\address[appliedmath]{Department of Applied Mathematics, University of Waterloo, Waterloo, Ontario N2L 3G1, Canada}
\address[hmi]{Hexagon Manufacturing Intelligence Canada, Oakville, Ontario L6L 0G4, Canada}
\address[mechanicalUW]{Department of Mechanical and Mechatronics Engineering,  University of Waterloo, Ontario N2L 3G1, Canada}
\address[yorktesserea]{Department of Civil Engineering, York University, Toronto, Ontario M3J 1P3, Canada}
\tnotetext[t1]{Research funded by a Discovery grant from the Natural Sciences and Engineering Research Council of Canada.}

\date{Received: date / Accepted: date}

\begin{abstract}
A 
method for quasistatic cohesive fracture is introduced that uses an alternating direction method of multipliers (ADMM) to implement an energy approach to cohesive fracture. 
The ADMM algorithm minimizes a non-smooth, non-convex potential functional 
at each strain increment to predict the evolution of a cohesive-elastic system. The optimization problem bypasses the explicit stress criterion of force-based (Newtonian) methods, which interferes with Newton iterations impeding convergence. The model is extended with an extrapolation method that significantly reduces the computation time of the sequence of optimizations. The ADMM algorithm is experimentally shown to have nearly linear time complexity and fast iteration times, allowing it to simulate much larger problems than were previously feasible. The effectiveness, as well as the insensitivity of the algorithm to its numerical parameters is demonstrated through examples. It is shown that the Lagrange multiplier method of ADMM is more effective than earlier Nitsche and continuation methods for quasistatic problems. Close spaced minima are identified in complicated microstructures and their effect discussed.

\end{abstract}


\end{frontmatter}



\section{Introduction}
Original approaches to fracture modeling were based on linear elastic fracture mechanics (LEFM), later generalized to the development of a plasticity zone at the crack tip. Based on failure criteria employing the concepts of stress intensity factor and the J-integral, they were typically used in a post-processing step assuming idealized geometries and then incorporated into finite element analyses. In contrast, cohesive fracture assumes gradual decohesion ahead of the fracture tip. This seems to be close to the reality of most materials. It also offers regularization of several aspects of LEFM, which treats each material point as either cracked or uncracked. Following the pioneering work of Griffith \cite{griffith1921phenomena}, several authors \cite{FraMar98,miehe2007configurational,del2013diffuse} cast LEFM as an energy minimization problem. Attendant numerical solutions gave rise to a class of smeared crack methods that employ a phase field \cite{bourdin2000numerical,miehe2007configurational}. These were later extended to include cohesive fracture \cite{geelenetal2019,wu2018length}. In a parallel development,  discrete cohesive fracture in the sense of initiation and propagation of a strong displacement discontinuity 
was also formulated as a minimization problem of a non-differentiable, non-convex energy functional in  \cite{Lorentz2008subgradient, energymethod}. The concept of subdifferential and a continuation method were used, respectively. 
Applied to cohesive models which postulate no separation prior to a critical stress combination (termed ``initially rigid'' or ``extrinsic'' cohesive models), these works are distinct from finite element methods that use the same weak form as obtained from minimization of a smooth functional. The latter methods deal with crack initiation through an externally imposed stress criterion in the course of solving the attendant ODEs while in non-smooth optimization methods the activation criterion lies within the optimization method in the sense that an interface will open when it is energetically favorable to do so and not according to an externally imposed criterion. It was shown in \cite{timecont} that the force based (as opposed to energy based) approach of solving equilibrium ODEs leads to a discontinuity of the internal force vector in time (termed a ``time discontinuity'') which results in an ill-posed set of ODEs. It was further shown that the discontinuity is the root cause of some undesirable features of cohesive crack analysis \cite{sam2005obtaining}. Although remedies specifically meant to address the time discontinuity problem were offered \cite{timecont,sam2005obtaining,chen2019nodal,chen2021nodal} as well as others, addressing undesirable behaviors without specifically attributing them to time discontinuity\cite{song2013explicit}, it was shown in \cite{energymethod,hirmand2019robust,CMaccepted,contDG} that time discontinuity can be cured with no additional effort by the energy method applied to a non-differentiable, non-convex functional.

The phase-field method has the advantage that it tracks the crack path independent of the element boundaries, whereas the discrete method, as implemented so-far \cite{energymethod,blockCD,contDG,vavasis2020second}, presupposes the crack path along element interfaces. This draw-back is currently addressed  by randomized \cite{leon2014reduction} or special isoperimetric \cite{pinwheel} meshes, the latter coming with mathematical guarantees of (slow) convergence in 2D. The method, however, can be implemented with XFEM or any other embedded strong discontinuity method. On the other hand, the discrete model is able to represent crack opening and therefore fragmentation, as well as multiphysics problems in which flow within the crack is part of the problem. 

In this paper we focus on the discrete cohesive energy method, 
whose avoidance of an explicit fracture criterion 
precisely leads to the above mentioned property of ``time continuity'', i.e., a model that  generates a continuous force vector at the time of crack activation.
It also endows the model with the ability to perform implicit calculations, whereas in force-based (as opposed to energy-based) methods, the criterion interferes with the Newton iterations of an implicit scheme thus impeding convergence. 

``Initially rigid'' or ``extrinsic'' cohesive models were initially proposed for explicit dynamics \cite{OrtizPan}. 
They are applied to quasistatic problems less commonly in the literature because of the difficulty that Newton's method does not converge on functions nondifferentiable at the solution.  Convergence of Newton's method has been studied in the context of nondifferentiable functions arising in fracture simulations using peridynamics  \cite{ni2019static}.  Other authors, for example \cite{liu2019quasi}, recognizing the issue, opted for an explicit solution of the quasistatic problem or circumvented difficulties of an implicit solution by either keeping the cracks fixed during Newton's iterations \cite{secchi2004cohesive,schrefler2006adaptive} or employing an ``elastic-cohesive'' (``intrinsic'') model \cite{jin2002finite}. 


The ability of the energy method to perform effective implicit calculations for quasistatic problems was confirmed by several implicit optimization implementations \cite{energymethod,contDG,vavasis2020second}, which, however, provided slow convergence.  Recent work 
\cite{blockCD} resulted in an 
optimization 
method, both implicit and explicit, 
based on a discontinuous Galerkin formulation with Nitsche flux
and a block coordinate descent (BCD) algorithm, which performs well in explicit dynamics but is slower in implicit calculations, especially for quasistatic problems. 
This paper 
develops a new iterative method using ADMM (Alternating Direction Method of Multipliers) for quasistatic fracture that is 2-3 orders of magnitude faster for the set of problems considered. The new method uses ADMM to optimize the non-convex, non-smooth, linearly constrained potential function sequentially for each level of load.

In Section \ref{energyapproach}, an overview of the cohesive-elastic potential function is given. Spatial discretization is obtained in Section \ref{spatial}. In Section \ref{ADMM}, the ADMM algorithm is applied to this problem and primal and dual convergence criteria are introduced. 
A new extrapolation technique, introduced in Section \ref{extrapolation}, is later (Section \ref{test}) shown to reduce the run-time by up to an order of magnitude by linearly estimating the next state and using this point as the optimization initial coordinate for the following quasistatic loading step.  This technique is believed to have potential use in other quasistatic models as well. An example problem in Section \ref{example} shows the dramatic improvement in run time of the proposed method compared to other methods.
The algorithm is dependent on four discretization and optimization parameters. To analyze the sensitivity of the results to these parameters, in Section \ref{test} several studies are performed on a test problem with adjusted parameter values. The time complexity of the ADMM algorithm is also experimentally tested. Nearly linear scaling with problem size and very fast iteration times allow the ADMM algorithm to solve large static fracture problems quickly.

\section{The energy approach and its spatial discretization}
\label{energyapproach}

We consider a class of models 
composed of a bulk linear-elastic material in which cohesive surfaces with non-zero opening form at a sufficient level of stress. This property of the cohesive cracks, i.e., an initially-rigid behaviour, is responsible for the singularity of the potential function 
and for the fact that implicit quasistatic solutions are less commonly encountered in the literature \cite{areias2008quasi,giovanardi2019fullycoupled,ferte20163d,peruzzo2019dynamics}. 
The system considered is a connected shape $\Omega$ in two dimensions (the logic generalizes to three dimensions) with a predefined set of interfaces $\Gamma_{d}$, taken to be finite element boundaries (see next section), interlacing the shape. 
The state of the system is given by the 
displacement field $\boldsymbol{u(x)}$, which admits discontinuities on $\Gamma_d$. The interface openings $
[\![\boldsymbol{u}(\boldsymbol{x})]\!]$ 
are defined as the jump discontinuity of $\boldsymbol{u(x)}$ along $\Gamma_d$,  
\begin{eqnarray}
[\![\boldsymbol{u}]\!]&=&\boldsymbol{u}^{+}-\boldsymbol{u}^{-}\\
\boldsymbol{u}^{ \pm}&=&\lim _{\varepsilon \rightarrow 0^{+}} \boldsymbol{u}\left(\boldsymbol{x} \pm \varepsilon \boldsymbol {n}_{d}\right),
\end{eqnarray}
and $\boldsymbol {n}_{d}$ is the unit normal on $\Gamma_d$, whose orientation is defined arbitrarily.  

\begin{figure}[ht]
\centering
\begin{subfigure}{0.47\textwidth}
    \centering
    \includegraphics[width=0.9\linewidth]{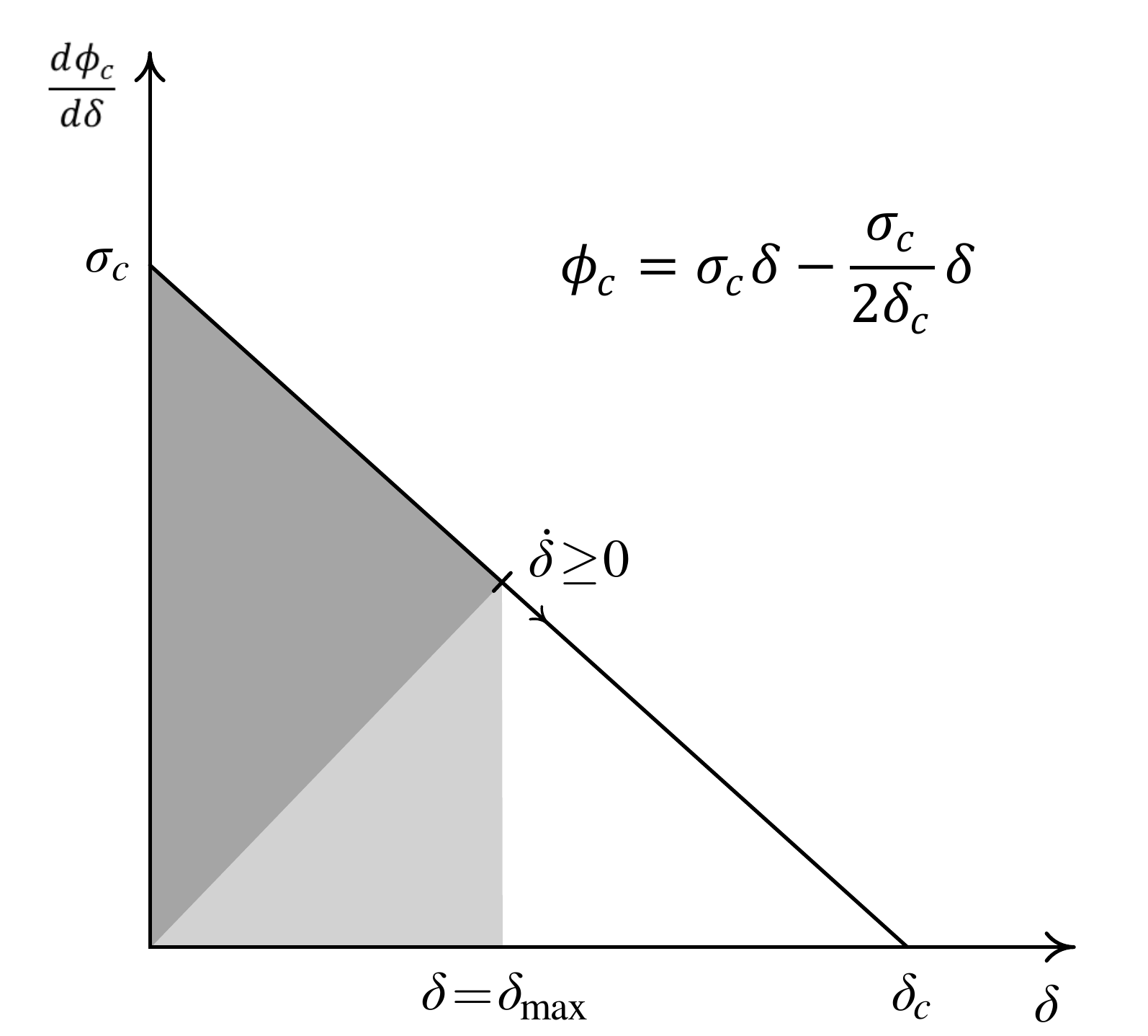}
    \caption{loading}
\end{subfigure}
\begin{subfigure}{0.52\textwidth}
    \centering
    \includegraphics[width=0.9\linewidth]{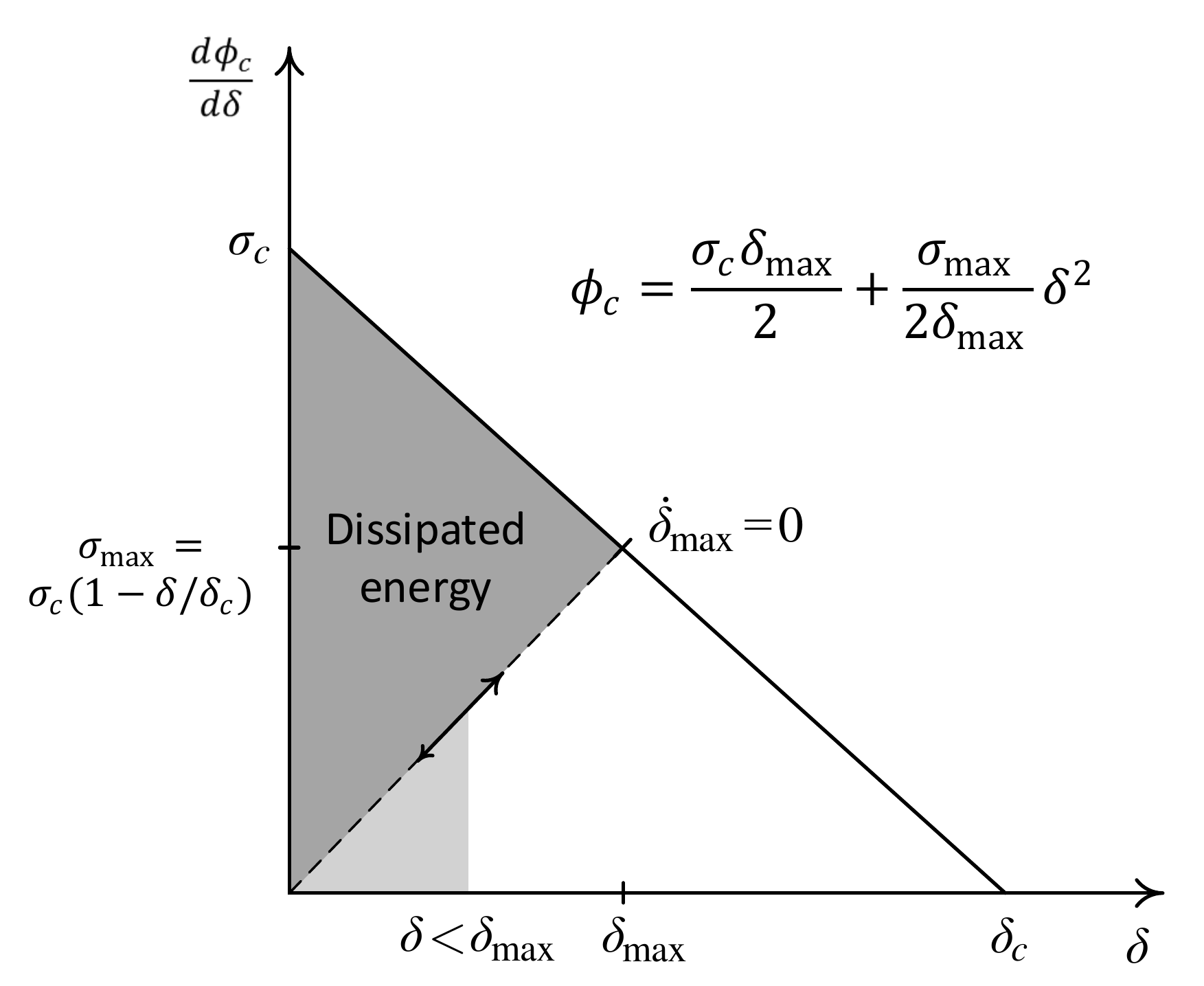}
    \caption{unloading/reloading}
\end{subfigure}
\caption{Initially rigid linearly decaying traction - separation model. $\delta = \sqrt{\delta_n^2 + \beta^2 \delta_s^2}$ is an effective opening displacement with ``mixity parameter'' $\beta$ as introduced in \cite{OrtizPan}; $\delta_{c}$ is the opening magnitude at complete failure, and $\sigma_{c}$ the critical stress that must be reached for failure initiation. The cohesive potential is irreversible: a) during initial loading (when $\delta > \delta_\text{max}$), some energy is dissipated; b) when unloaded and re-loaded to a displacement less than $\delta_\text{max}$, stress increases proportional to $\delta$. In each case, $\phi_c$ is equal to the total shaded area. The dissipated portion of the energy is shaded darker. }
\label{fig:cohesivemodel}
\end{figure}

\begin{figure}[ht]
\centering
    \centering
    \includegraphics[width=0.6\linewidth]{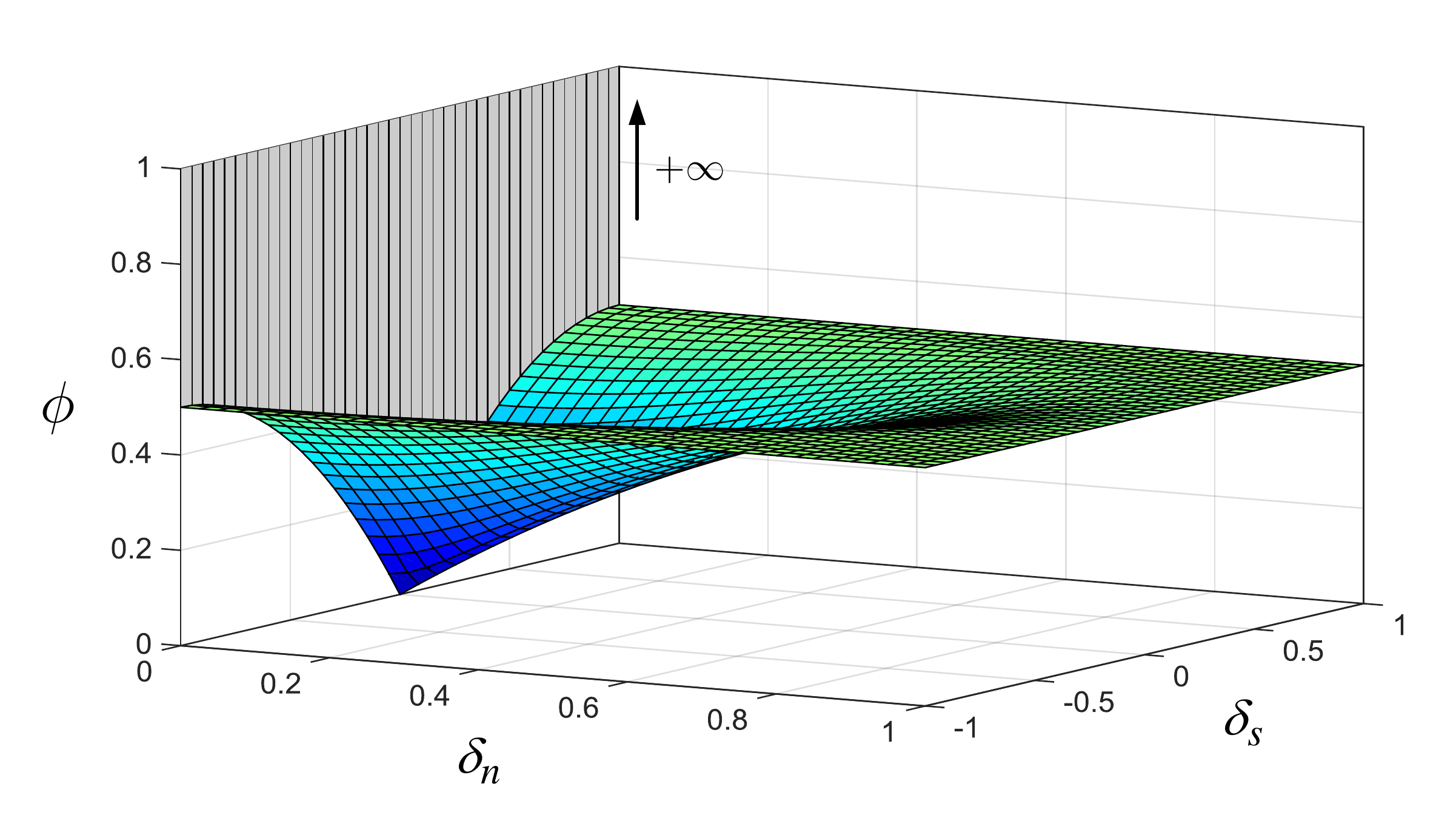}
    
\caption{Variation of the cohesive potential against $(\delta_n, \delta_s)$ with $\sigma_c=1$, $\delta_c=1$, and $\beta=1$, capped off at $\delta_n<0$.  It is assumed that the interface undergoes loading and
no unloading occurs.}
\label{fig:totalpotential}
\end{figure}

The potential energy of the system assuming no body forces is defined 
as
\begin{equation}
\label{potential}
\pi(\boldsymbol{u})=\int_{\Omega \backslash \Gamma_{d}} \psi(\boldsymbol{\varepsilon}(\boldsymbol{u})) \mathrm{d} V+\int_{\Gamma_{d}}\phi(\boldsymbol{\delta})
\mathrm{d} S 
- \int_{\partial_t \Omega} \boldsymbol{u} \cdot \bar{\boldsymbol{t}} \mathrm{d} S\\
\end{equation}
where 
\begin{eqnarray}
    \boldsymbol{\delta}&=&[\![\boldsymbol{u}(\boldsymbol{x})]\!],\\
    \phi(\boldsymbol{\delta})&=&\phi_c (\delta) + I_{\mathbb{R}^+} (\delta_n),
\end{eqnarray}
$\phi_c (\delta)$
is  the cohesive potential function, an example of which for a linearly descending cohesive interface model is defined in Fig. \ref{fig:cohesivemodel}; $\delta$ is the magnitude of $\boldsymbol{\delta}=(\delta_n,\delta_s)$, also defined in  Fig. \ref{fig:cohesivemodel}; 
$\delta_n$ is the normal component of the interface opening vector in the direction of $\boldsymbol {n}_{d}$ and $\delta_s$ the component of the same vector in the direction parallel to $\Gamma_d$;  $I_{\mathbb{R}^+} (\delta_n)$ is an indicator function (in the sense of convex analysis \cite{clarke1990optimization}) enforcing that no interpenetration occurs along $\Gamma_d$ by assigning infinite energy to the non-feasible configuration $\delta_n<0$; $\bar{\boldsymbol{t}}$ is a vector of tractions applied to part of the body's boundary denoted $\partial_t \Omega$. 
We note that both $\phi_c(\delta)$ and $I_{\mathbb{R}^+} (\delta_n)$ are not differentiable at crack initiation corresponding to $\delta=\delta_\text{max}=0$; see Fig. \ref{fig:totalpotential}. 
$\psi(\boldsymbol{\varepsilon}(\boldsymbol{u}))$ is the linear-elastic energy associated with displacement $\boldsymbol{u}$. 
The material is assumed to behave in a linear-elastic manner so the local elastic energy density obeys:
\begin{equation}
\psi(\boldsymbol{\varepsilon}(\boldsymbol{u})) = 1/2 \boldsymbol{\varepsilon}(\boldsymbol{u}) : D : \boldsymbol{\varepsilon}(\boldsymbol{u}),
\end{equation}
where $D$ is the elasticity constitutive tensor and $\boldsymbol{\varepsilon}(\boldsymbol{u})$ is the symmetric part of the displacement gradient.
In a typical problem a subset of $\boldsymbol{u}$ will be set by the Dirichlet boundary conditions.


\label{spatial}
The optimization problem is spatially discretized in 2D plane stress/plane strain applications by a triangular finite element mesh  with 
multiple copies of nodes so that no element shares a node with a neighbouring element. Cohesive interface elements are positioned between elements at each element boundary at the start of the simulation acting as potential sites of crack initiation and propagation (i.e., $\Gamma_d$ is taken to be the union of all element edges in the mesh). 
The deformation state is parameterized by $\boldsymbol{u}_h$, a vector containing the global displacements of the nodes
and an interface opening displacement $\boldsymbol{\delta}_i$ per Gauss point $i$ of each interface surface, computed as an additional variable and constrained to equal $[\![\boldsymbol{u}_i]\!]$ at the Gauss point; $[\![\boldsymbol{u}_i]\!]$ and  $\boldsymbol{\delta}_i$ can be computed in either the global or the local coordinate system of the interface (a normal and a parallel component) at the Gauss point. The initial surface orientation is used to define the normal direction and is not updated during the simulation (for use in problems with significant rotation this should be updated). In either case, the operation $\boldsymbol{u}_h \rightarrow [\![\boldsymbol{u}_i]\!]$ is a linear transformation composed of interpolation, vector difference, and rotation. The operator $A_i$ is used to denote the linear mapping at Gauss point $i$: 

\begin{equation}
\label{linear}
[\![\boldsymbol{u}_i]\!] = A_i\boldsymbol{u}_h.
\end{equation}

Let $\boldsymbol{\delta}_h$ define the set of interface opening vectors at all interface Gauss points. Let $K$ be the global stiffness matrix of the triangles (i.e., the ``bulk'' finite elements) and consider a Gauss quadrature is employed for numerical integration on $\Gamma_d$. The total energy of the discretized system in the absence of body forces and surface tractions is
 \begin{equation}
     \pi({{\boldsymbol u}_h, \boldsymbol\delta}_h)= E(\boldsymbol{u}_h) +
 C(\boldsymbol{\delta}_h),   
 \end{equation}
where
\begin{equation}
\label{elasticen}
E(\boldsymbol{u}_h) = (1/2) \boldsymbol{u}_h^TK\boldsymbol{u}_h
\end{equation}
is the total elastic energy of the discrete system and 
\begin{equation}
 C(\boldsymbol{\delta}_h) = \sum_{i \in IntfPoints}  a_i\phi_c (\delta_i)
 \end{equation}
is the cohesive energy, with $a_i$ the effective (weighted) area of interface point $i$ and $\delta_i = \sqrt{\delta_{n_i}^2 + \beta^2 \delta_{s_i}^2}$, formed from the 2-vector $\boldsymbol{\delta}_{h_i} = (\delta_{n_i}, \delta_{s_i})$ within $\boldsymbol{\delta}_h$ corresponding to interface point $i$. In what follows the subscript $h$ pertaining to spatial discretization will be dropped, i.e., vectors $\boldsymbol{u}, \boldsymbol{\delta}$ will be used for $\boldsymbol{u}_h, \boldsymbol{\delta}_h$, respectively, not to be confused with the 2-vectors $\boldsymbol{u}$, $\boldsymbol{\delta}$ of the continuum formulation. Replacing equation (\ref{linear}), $A\boldsymbol{u}$ will denote the displacement jumps obtained from ${\boldsymbol u}_h$ at interface Gauss points, collectively.

\section{Quasi-static minimization and ADMM algorithm}
\label{ADMM}

A quasi-static process is a process that occurs slowly enough so as to always be in equilibrium. An equilibrium state 
is one of locally minimal energy. It is therefore natural for a  time discrete mathematical model of a quasi-static process to be posed as a sequence of energy minimization problems subject to constraints defined by the changing boundary conditions. 
The problem statement can therefore be defined as:
\begin{equation}
\label{eq:mini-problem}
({\boldsymbol u}, \boldsymbol{\delta})^{k} = {\rm argmin}(\pi({\boldsymbol u}, \boldsymbol{\delta})) \quad \text{s.t.} \quad \boldsymbol{u}_{BC} =\bar{\boldsymbol{u}}_{BC^k}, \>\>\boldsymbol{\delta} = A{\boldsymbol u}
\end{equation}
where ${\boldsymbol u}, \boldsymbol\delta$ are nodal displacements and interface point openings; $\boldsymbol{u}^k$ is the solution at load step $k$, $\boldsymbol{u}_{BC}$ is the subset of $\boldsymbol{u}$, at which the Dirichlet constraints are applied, $\bar{\boldsymbol{u}}_{BC^k}$ is the value of the boundary conditions at load step $k$. 
It is assumed that the system is under such conditions that it will equilibrate at the nearest local minimum.
Because the cohesive energy function is non-convex, the entire functional is non-convex, which means that there will in general be multiple local minimizers.  The optimization method proposed here will converge to one of these. The computed minimizer  is  close to the equilibrium from the previous load step in most cases.  This is physically plausible and, as reported, 
 the energy minimization computations using different methods \cite{contDG,blockCD,vavasis2020second} match physical experiments. 
However, the underlying physical system may be fundamentally unstable in the sense that a small perturbation to the geometry or parameters leads to a significantly different solution path. One such occurence of importance to multiscale analysis is identified in the test problem of Section \ref{test}. This behavior can be detected by an optimization method like the one proposed here only by sampling the space of geometries and parameters with multiple runs.  


With algorithmic speed and simplicity as 
priorities, ADMM was chosen as the optimization algorithm. In \cite{admmconv} it is shown that ADMM will converge for a class of non-smooth, non-convex functions quite similar to the potential functions considered here. As an iterative method, ADMM simplifies the global optimization problem into one that is an unconstrained quadratic optimization problem and a set of low dimensional optimization problems that can be solved analytically. 
This aspect is similar to \cite{blockCD}, except that in that method  the explicitly parametrized Lagrange multipliers used here (see below) are replaced with the Nitsche flux, which has advantages in a dynamic setting but is not as effective for quasistatic problems. 
A significant difference between optimization for dynamic problems and optimization for quasistatic problems is that quasistatic optimization loses the regularization provided by the mass matrix. 
The continuity of tractions at interface points at the time of crack activation 
\cite{timecont} remains relevant to ensure a smooth development of forces with evolution of the deformation throughout the solution trajectory. This is achieved with the use of the energy method as explained in \cite{energymethod}. Quasistatic iterative solutions with a stress-based (Newtonian) formulation, i.e., one that relies on extrinsic fracture criteria to initiate and propagate fractures, 
are indeed very difficult to converge \cite{hirmand2019nondifferentiable}.

In the following, the ADMM algorithm is presented and primal and dual convergence checks \cite{boyd2011distributed}, computed using the infinity norm on residuals with units of pressure, are introduced. First, the augmented Lagrangian for the problem is defined:
\begin{equation}
L(\boldsymbol{u},\boldsymbol{\delta};\boldsymbol{y}) =  E(\boldsymbol{u}) + C(\boldsymbol{\delta}) + \boldsymbol{y}^T(A\boldsymbol{u}-\boldsymbol{\delta}) + (\rho/2)||A\boldsymbol{u} - \boldsymbol{\delta}||^2,
\end{equation}
where $\boldsymbol{y}$ is the Lagrange multiplier vector of the same length as $\boldsymbol\delta$ (i.e., one Lagrange multiplier $\boldsymbol{y}_i$ per interface Gauss point $i$) and $\rho$ is an optimization parameter (penalty). The ADMM algorithm \cite{boyd2011distributed} is: 

\vspace{10pt}
\begin{algorithmic} 
\WHILE{$!converged$}
\STATE $\boldsymbol{u}^{k+1} = min_{\boldsymbol{u}}[E(\boldsymbol{u}) + (\boldsymbol{y}^{k})^T(A\boldsymbol{u}) + (\rho/2) \>||A\boldsymbol{u} - \boldsymbol{\delta}^{k}||^2]$
  \STATE $\boldsymbol{\delta}^{k+1} = min_{\boldsymbol{\delta}}[C(\boldsymbol{\delta}) + (\boldsymbol{y}^{k})^T(-\boldsymbol{\delta}) + (\rho/2)\>||A\boldsymbol{u}^{k+1} - \boldsymbol{\delta}||^2]$
  \STATE $\boldsymbol{y}^{k+1} = \boldsymbol{y}^{k} + \rho \> (A\boldsymbol{u}^{k+1} - \boldsymbol{\delta}^{k+1})$
\ENDWHILE
\end{algorithmic}
\vspace{10pt}

This statement of the problem 
defines the displacement variable and the fracture opening variable at an interface point as separate values and then specifies their relationship using the linear constraint $A\boldsymbol{u} = \boldsymbol{\delta}$. The benefit of this approach is that the separate optimizations in $\boldsymbol{u}$ and $\boldsymbol{\delta}$ are much simpler. For the deformation minimization, this is exactly the same as for a normal elastic problem and can be solved by 
solving a linear system: with $E({\boldsymbol u})$ given by (\ref{elasticen}), the $\boldsymbol u$ minimization
expands to 
\begin{equation}
{\boldsymbol u}^{k+1} = min_{\boldsymbol u}[(1/2){\boldsymbol u}^T(K+\rho A^TA) {\boldsymbol u} + {\boldsymbol u}^T A^T ({\boldsymbol y}^k - \rho {\boldsymbol{\delta}}^k) ],
\end{equation}
requiring
\begin{equation}
 (K + \rho A^TA)u^{k+1} + A^T(\boldsymbol{y}^{k} - \rho \boldsymbol{\delta}^{k}) = 0   
\end{equation}
or
\begin{equation}
{\boldsymbol u}^{k+1} =  - (K + \rho A^TA)^{-1} A^T({\boldsymbol y}^{k} - \rho \boldsymbol{\delta}^k),
\end{equation}
assuming that $(K + \rho A^TA)$ is positive definite and therefore invertible, which is the case for any positive, nonzero penalty parameter $\rho$. 

The crack opening minimization is written as
\begin{eqnarray}
\label{eq:localmin1}
{\boldsymbol{\delta}}^{k+1} &=& min_{\boldsymbol{\delta}}[C(\boldsymbol{\delta}) - ({\boldsymbol y}^{k})^T\boldsymbol{\delta} + (\rho/2) \>||A{\boldsymbol u}^{k+1} - \boldsymbol{\delta}||^2].
\end{eqnarray}

Due to separability of the unknown opening displacements $\boldsymbol{\delta}_i$, this becomes a set of many independent minimization problems:
\begin{equation}
\label{eq:localmin}
min_{\boldsymbol{\delta}_i} \left[ a_i \phi_c(\delta_i) - (\boldsymbol{y}^{k}_i)^T\boldsymbol{\delta}_i + (\rho/2) \>||A_i{\boldsymbol u}_i^{k+1} - \boldsymbol{\delta}_i||^2  \right] \quad \forall i \in IntfPoints.
\end{equation}

\begin{figure}[ht]
\centering
\begin{subfigure}{0.45\textwidth}
    \centering
    \includegraphics[width=0.9\linewidth]{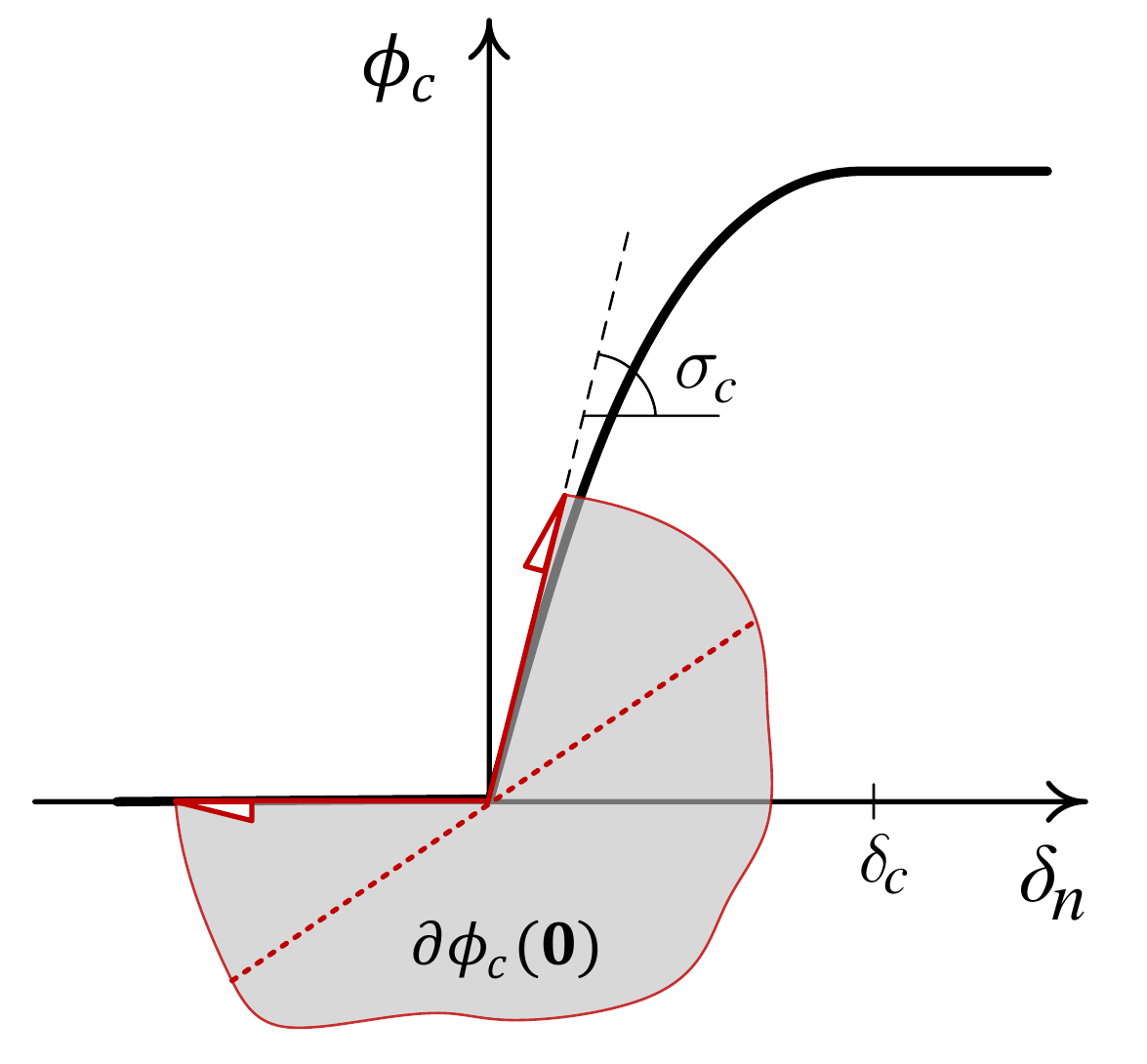}
    \caption{}
    \label{fig:cohesivemodel-a}
\end{subfigure}
\begin{subfigure}{0.45\textwidth}
    \centering
    \includegraphics[width=1.0\linewidth]{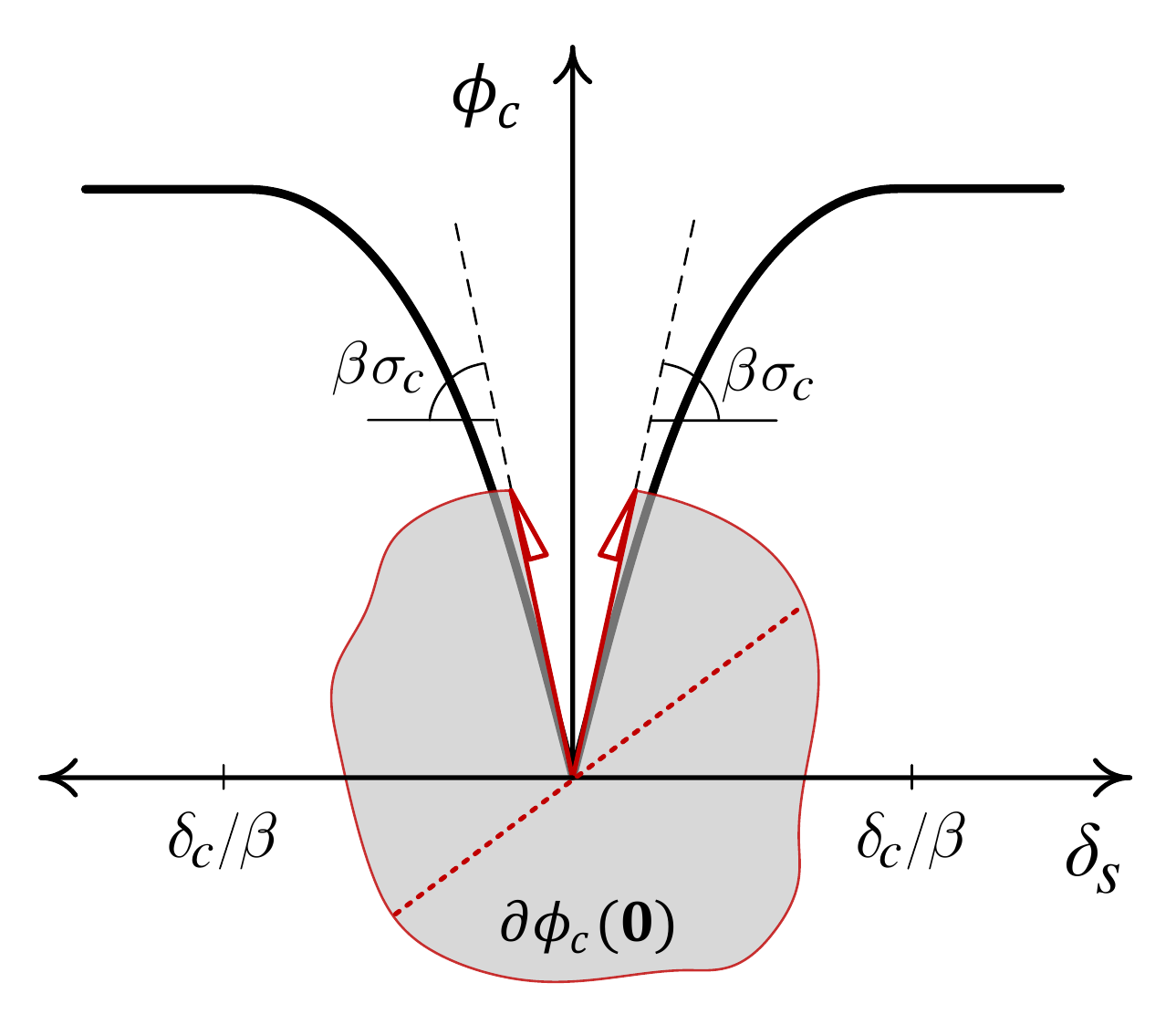}
    \caption{}
    \label{fig:cohesivemodel-b}
\end{subfigure}
\caption{Graphical interpretation of $\phi_c$ at the origin; a) in the plane $\delta_s = 0$ and b) in the plane $\delta_n = 0$.}
\label{fig:generaldiff-phi}
\end{figure}

\begin{figure}[ht]
\centering
    \centering
    \includegraphics[width=0.4\linewidth]{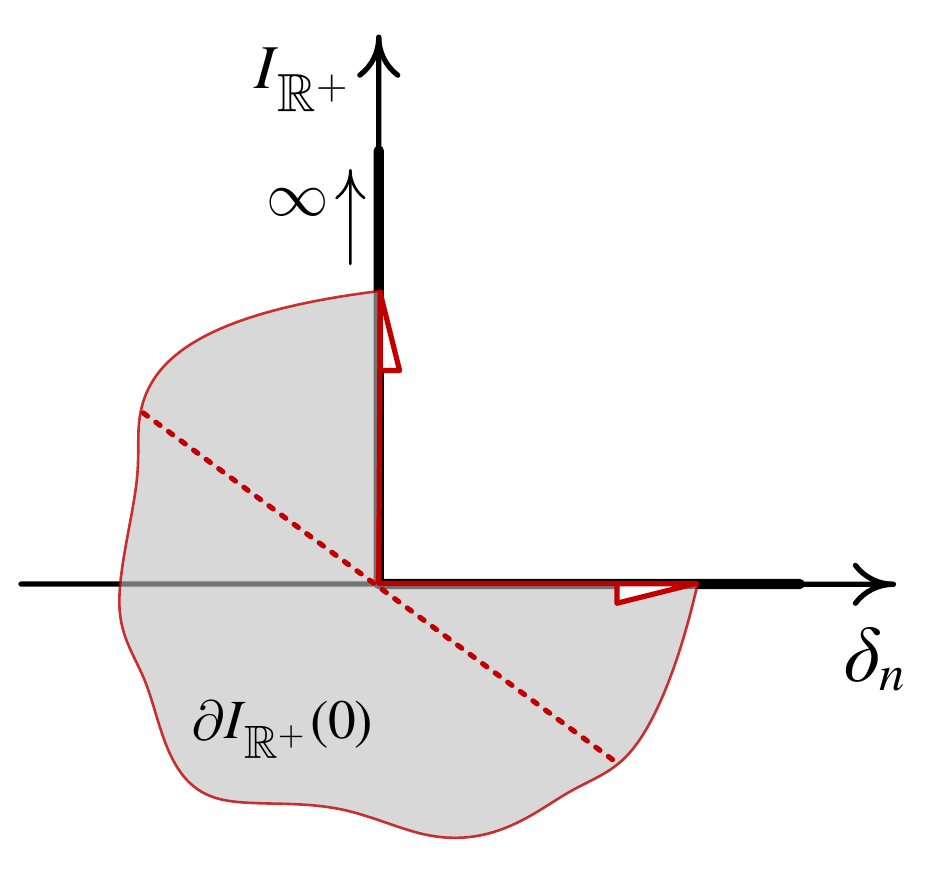}
\caption{Graphical interpretation of $\partial I_{\mathbb{R}^+}$.}
\label{fig:generaldiff-indicator}
\end{figure}

This minimization entails the derivative of $\phi_c(\delta_i)$ w.r.t. $\boldsymbol\delta$. Because $\phi_c(\delta_i)$ is not globally differentiable, a straightforward solution of the minimization problems \eqref{eq:localmin} is out of reach. A computational algorithm 
is employed which exploits a generalized differential calculus to solve the non-differentiable problem \eqref{eq:localmin}. The generalized gradient of a function $f$ at a non-differentiable point $\boldsymbol{x}$, which is denoted $\partial f(\boldsymbol{x})$, is defined as the set of slopes that are less steep than the slope of any directional derivative of $f$ at $\boldsymbol{x}$ considering all admissible directions. Graphical interpretations of $\partial \phi_c$ and $\partial I_{\mathbb{R}^+}$ at the origin are schematically shown in Fig. \ref{fig:generaldiff-phi} and \ref{fig:generaldiff-indicator}. At differentiable points, the generalized gradient reduces to the singleton set $\{\nabla f(\boldsymbol{x})\}$. Just as in differential calculus, where derivatives are rarely computed from the definition, one appeals to a body of theory and to certain rules that characterize generalized gradients, see \cite{clarke1990optimization}. Using generalized differential calculus, one writes: $\boldsymbol{\delta}^\star$ is a solution of \eqref{eq:localmin} if and only if

\begin{equation}
    \label{eq:generalzied-diff-cond}
    \mathbf{0} \in \partial \phi_c(\boldsymbol{\delta}^\star_i) + \partial I_{\mathbb{R}^+}(\delta^\star_{n,i}) + \{\rho \boldsymbol{\delta}^\star_i - \boldsymbol{p}_i\},
\end{equation}
where the summation is performed  in a \emph{Minkowski} sense (see \cite{clarke1990optimization} for a definition of Minkowski sum) and $\boldsymbol{p}_i$ is a traction computed from the known deformation and Lagrange multipliers at interface point $i$,
\begin{equation}
    \boldsymbol{p}_i= \boldsymbol{y}_i + \rho \boldsymbol{A}_i \boldsymbol{u}_i.
\end{equation}
Equation \eqref{eq:generalzied-diff-cond} can be used to readily determine if the minimizer occurs at a non-differentiable point. If not, solution is sought in a domain where the objective is differentiable, which can be solved analytically for $\beta = 1$ (or numerically otherwise). We refer the reader to \cite{blockCD} for a detailed derivation of these criteria for brevity of presentation. The resulting algorithm is:

\vspace{10pt}
\begin{algorithmic}
\STATE Compute $p_i=\sqrt{\max(p_{n,i},0)^2+\|\mathbf{p}_{s,i}/\beta\|^2}$
\STATE \textbf{if} $p_i\le\sigma_c$  and $\delta_{\text{max},i}=0$, \textbf{then}
\STATE \hspace{5mm} pre-activation state with 
$\boldsymbol{\delta}_i=\mathbf{0}$
\STATE \textbf{else}
\STATE \hspace{5mm}
if $p_{n,i}\le0$, set $\delta_{n,i}=0$ and solve  \eqref{eq:localmin} for $\boldsymbol{\delta}_{s,i}$ only,
\STATE \hspace{5mm} else, declare $\delta_{n,i}>0$ and solve \eqref{eq:localmin} for $\boldsymbol{\delta}_i$ 
\STATE \textbf{end}
\end{algorithmic}
\vspace{10pt}
.

The remaining information needed to complete the model definition is the choice of the augmented Lagrange parameter, $\rho$, and that of convergence criteria. For each individual crack opening energy minimization to have a unique solution, the problem must be strongly convex. This can be achieved  \cite{blockCD} if
\begin{equation}
\label{eqn:rhoCriteria}
\rho >  {a_i} \frac {\sigma_c} {\delta_c} 
 \; \forall i. 
 \end{equation}
Larger values of $\rho$ effectively limit the distance travelled in each optimization iteration. This parameter can be used to make the optimization more conservative and not escape local minima as easily. For the following tests a value was chosen experimentally, where the trade-off between run-time and accuracy was weighed (larger values typically resulted in more consistent solutions but required more optimization iterations).

The convergence criteria used are a modification of those suggested in \cite{admmcrit}. The optimality conditions are primal and dual feasibility; primal feasibility can be computed directly while dual feasibility can be estimated using the current and previous state vectors. They specify that the linear constraints are satisfied and that the projection of the generalized gradient onto the feasible plane contains the zero vector. In accordance, the following residuals are formed:
\begin{align}
\boldsymbol {r}^{k+1}&=A \boldsymbol{u}^{k+1} - \boldsymbol {\delta}^{k+1}
\\
\boldsymbol {s}^{k+1}&=\rho A^{T}\left(\boldsymbol{\delta}^{k+1}-\boldsymbol {\delta}^{k}\right).
\end{align}

The primal residual, as defined above, has units of distance. By multiplying by $\rho$ and dividing each 
entry by the effective area of the corresponding interface Gauss point, the quantity is converted to a set of fictional pressures being applied by the augmented component of the Lagrangian to enforce the equality constraint. The convergence check used is the infinity norm of this vector of pressures. By taking the maximal value, this check can be considered a guarantee of local equilibrium. If, instead, the sum of pressures is used, then the check becomes increasingly more stringent with higher mesh resolution. Taking the average augmented pressure has a similar problem: creating a larger mesh with areas not under stress would reduce the importance of the areas locally out of equilibrium, making the convergence check less stringent: 
\begin{equation}
\tilde {\boldsymbol r}^{k+1}=\rho(A \boldsymbol {u}^{k+1} - \boldsymbol {\delta}^{k+1})\oslash \boldsymbol {a},
\end{equation}
where $\oslash \boldsymbol a$ denotes element-wise division by the effective interface areas. A similar transformation is applied to the dual residual, which originally has units of force. By dividing each 2-vector in $\boldsymbol s$ by the effective surface area of the corresponding 
interface point,
the vector of forces 
is converted to a vector of pressures. Each of these pressures gives an estimate of the local energy imbalance. This can be interpreted as an estimate of the potential gradient projected onto the constraint surface and normalized by the effective surface area. For the same reasons as above, the convergence check uses the infinity norm of this dual pressure residual:
\begin{equation}
\tilde {\boldsymbol s}^{k+1}=\rho A^{T}[(\boldsymbol {\delta}^{(k+1)}-\boldsymbol {\delta}^{(k)})\oslash \boldsymbol {a}].
\end{equation}
Then, in order to move to the next quasistatic step the following two conditions must be satisfied:
\begin{align}
||\tilde {\boldsymbol r}^{k+1} ||_{\infty} &< c_{\rm primal}
\label{primalPressure}
\\
||\tilde {\boldsymbol s}^{k+1} ||_{\infty} &< c_{\rm dual}.
\label{dualPressure}
\end{align}

The convergence criteria constants are defined based on the requirements of the problem. Due to their construction, they do not have to be updated for changes in $\rho$, mesh size, or number of elements. The maximum physical pressure experienced in a simulation is bounded above by $\sigma_c$, and the typical pressure of most areas of interest is usually close to this value. So requiring that $c_{dual} < \sigma_c/300$ and $c_{primal} < \sigma_c/300$ was found to give reasonably converged solutions. Because $\rho$ must be large to ensure small steps are taken and that the $\delta$ optimization is convex, the primal residual is typically an order of magnitude smaller than the dual residual. Usually, when this is observed in an ADMM optimization, $\rho$ would be decreased to achieve a more balanced convergence rate, but for the test problem in Section \ref{test}, the slight inefficiency is accepted.

\section{Extrapolation method}
\label{extrapolation}
An 
extrapolation technique was used to reduce computation time in subsequent optimizations. The default choice for the starting point for the next load step is the converged solution of the previous load step. This is probably the safest choice in general, however, it does not take advantage of the local smoothness often seen in the problem. Of the steps taken in the simulation, a majority do not experience a significant change in the number of open cracks. A reasonable assumption in this situation is that the system is locally almost linear, 
therefore the updated state in a new load step will be well approximated by a linear extrapolation of the previous converged solutions. Let $\boldsymbol{z}^k=[\boldsymbol{u}^k;\boldsymbol\delta^k;\boldsymbol{y}^k]$ denote the state of the optimization system at step $k$ and let $\tilde{\boldsymbol{z}}^{k+1}$ be the estimate for the next time step. Since the problem constraints are linear and both $\boldsymbol{z}^k$ and $\boldsymbol{z}^{k-1}$ are feasible states, any linear combination of them will also be feasible, i.e., it will have low primal residual.
Then the approximation being made is
\begin{equation}
\tilde{\boldsymbol{z}}^{k+1} = \boldsymbol{z}^k + (\boldsymbol{z}^k - \boldsymbol{z}^{k-1})
\end{equation}
and $\tilde{\boldsymbol{z}}^{k+1}$ is used as the initial estimate for the $k+1$ optimization step. This is done by setting the initial state to the extrapolated value, updating the loading conditions, and then iteratively solving for the actual updated state, $\boldsymbol{z}^{k+1}$. The decision of whether to use an extrapolation estimate is made based on the accuracy of the extrapolation in the previous step, as given by 
\begin{equation}
\label{eqn:ExtrapAcc}
||\boldsymbol{z}^{k} - \boldsymbol{z}^{k-1}|| / ||\boldsymbol{z}^{k} - \tilde{\boldsymbol{z}}^{k}||.
\end{equation}
If the current guess quality is better than a set threshold, for example
greater than 2, then an extrapolation estimate is used for the next step. Otherwise the current point, $\boldsymbol{z}^{k}$, is used as the initial state. 

In Section \ref{ParamSensitivity}, the impact of the extrapolation method on the solution and number of optimization iterations is investigated for an example problem. Generally, depending on the problem, extrapolation was observed to provide between zero improvement and a factor of 10 speedup. In some cases, slightly different results have been predicted when using the extrapolation estimate. This is believed to be caused by jumping into a neighbouring local minimum. \par
The extrapolation method is not specific to ADMM and could likely also be applied to other models, in which consecutive optimizations are being performed. Potential improvements are
the use of higher order extrapolation methods or learned heuristics; 
testing a small set of predicted starting points and picking the one with the lowest residual; and
checking the energy at intermediate points between $\tilde{\boldsymbol{z}}^{k+1}$ and $\boldsymbol{z}^{k}$ to ensure a local minimum is not being escaped.



\section{Example problem}
\label{example}

\begin{figure}
\begin{subfigure}{0.9\textwidth}
\begin{center}
\includegraphics[max width=0.8\linewidth]{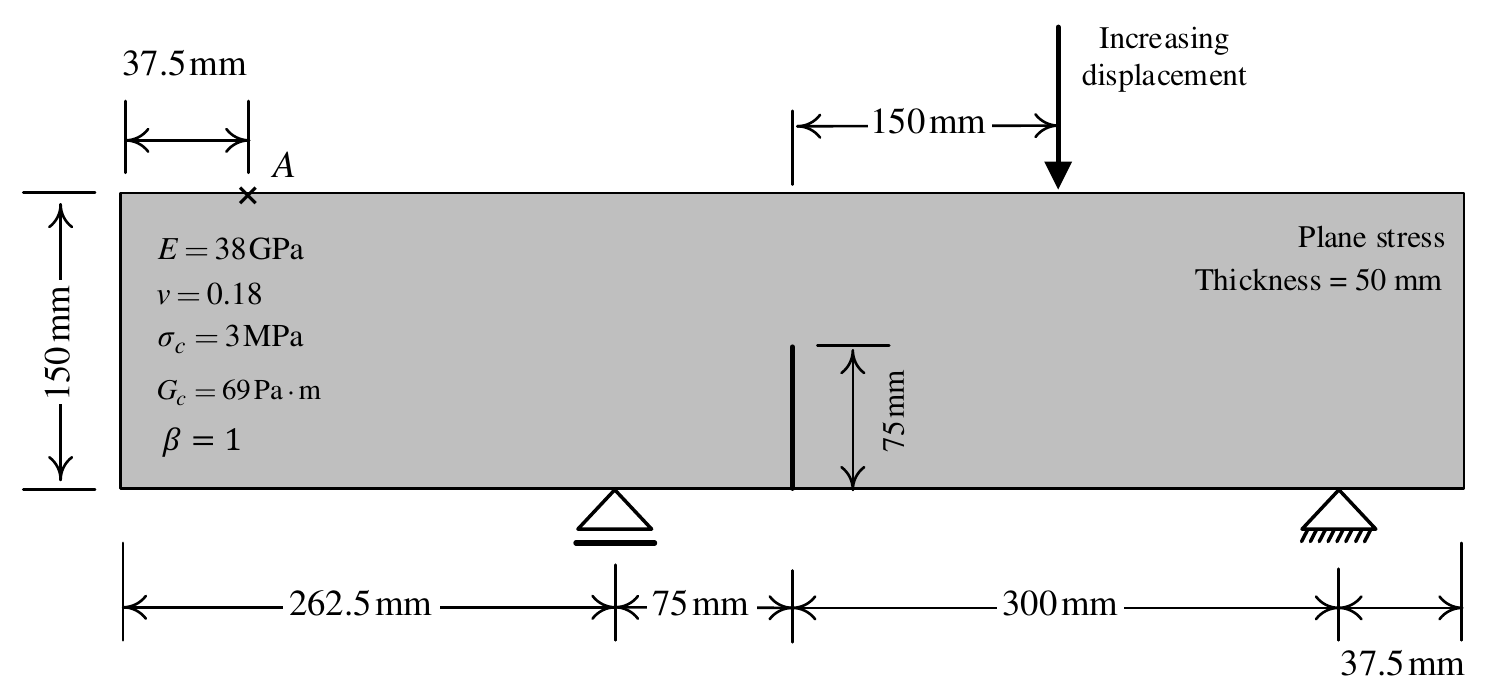}
\caption{}
\end{center}
\end{subfigure}

\begin{subfigure}{0.9\textwidth}
\begin{center}
\includegraphics[max width=0.75\linewidth]{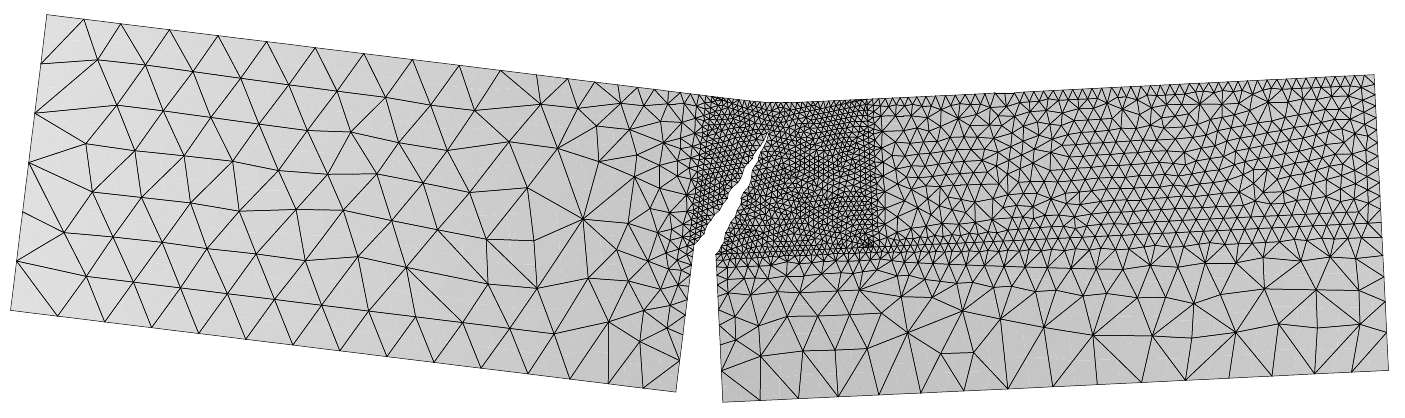}
\caption{}
\end{center}
\end{subfigure}
\caption{Single-edge notched beam; a) problem setting and relevant mechanical properties and, b) deformed geometry (magnified by a factor of 100) of the single edge notched beam with 0.075mm applied vertical displacement.}
\label{fig:galvez}
\end{figure}

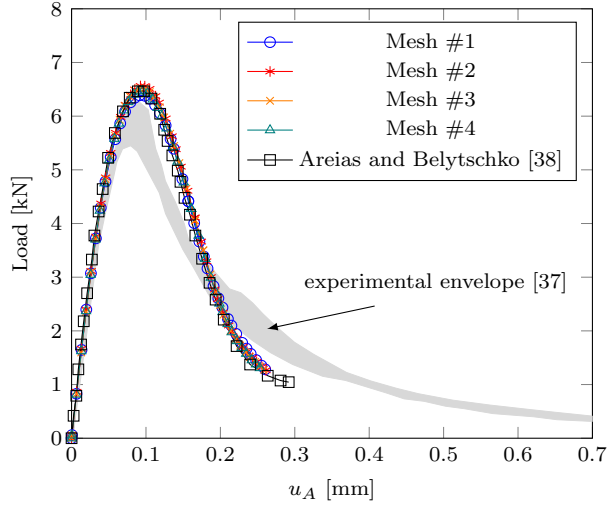
\begin{figure}[ht]
\centering
    \begin{tikzpicture} 
    \begin{axis}[
        xlabel={$u_A$ [mm]},
        ylabel={Load [kN]},
        xmin=0, xmax=0.7,
        ymin=0, ymax=8,
        xtick={0,0.1,0.2,0.3,0.4,0.5,0.6,0.7},
        ytick={0,1,2,3,4,5,6,7,8},
        legend pos=north east,
        /tikz/font=\footnotesize,
        legend style={font=\footnotesize}
    ]
        
     \addplot[forget plot] graphics
       [xmin=0,xmax=0.7,ymin=0,ymax=6.31]
       {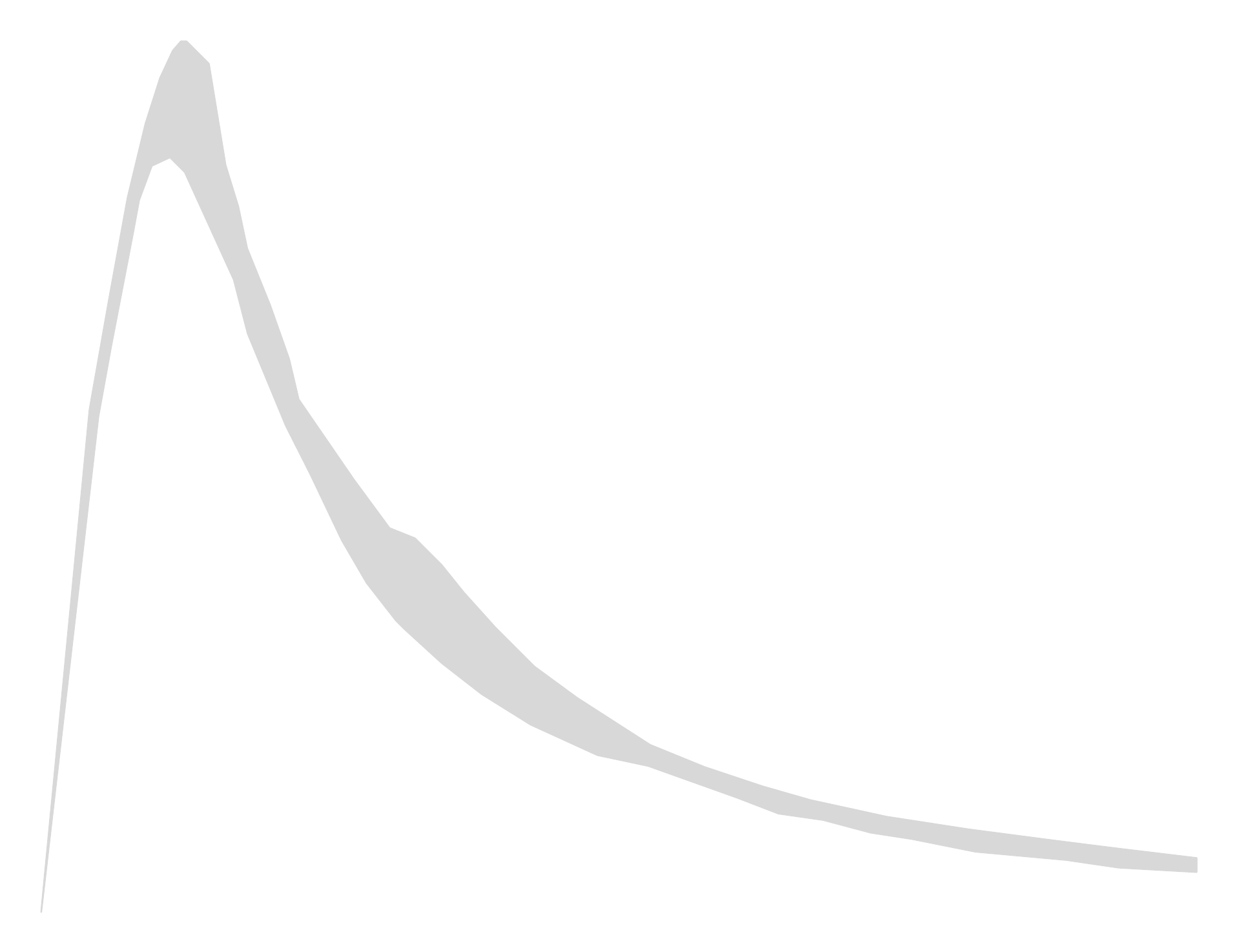};    
    
    \addplot[color=blue,  mark={o}]
        table {mesh1.dat};

    \addplot[color=red, mark={asterisk}]
        table {mesh2.dat};

    \addplot[color=orange, mark={x}]
      table {mesh3.dat};
        
    \addplot[color=teal, mark=triangle]
        table {mesh4.dat};

    \addplot[color=black, mark=square]
        table {areias.dat};

    \legend{Mesh $\#$1, Mesh $\#$2, Mesh $\#$3, Mesh $\#$4, Areias and Belytschko \cite{areias2005analysis}}
    
    \node (A) at (0.42, 2.5) {};
    \node (B) at (0.25, 2.0) {};
    \node at (0.49, 2.9) {experimental envelope \cite{galvez1998mixed}};
    \draw [-latex] (A) -- (B);
    
    \end{axis}
    \end{tikzpicture}
\caption{Load vs. displacement of point A of the single-edge notched beam obtained using different meshes}
\label{fig:galvez-load-disp-meshes}
\end{figure}

We compare the performance of the ADMM algorithm with algorithms in our previous work including block-CD \cite{blockCD}, continuation \cite{contDG}, and interior point \cite{vavasis2020second} algorithms. The problem considered for this purpose is the mixed-mode failure of the single-edge notched beam of Galvez et al. \cite{galvez1998mixed} which is summarized in Figure \ref{fig:galvez}. The FE mesh used was the same as the finer mesh in \cite{contDG} consisting of 3508 6-noded triangular elements. All algorithms yield seemingly identical results which are shown in Figure \ref{fig:galvez-load-disp-meshes} along with a comparison with experimental observations reported in \cite{galvez1998mixed} and numerical simulations obtained in \cite{areias2005analysis} using a 3D XFEM model. The run-times of different algorithms are summarized in Table \ref{tab:comparison}. 
All methods were implemented in our in-house code and convergence was judged by considering similar error tolerances. As opposed to the ADMM and block-CD algorithms, the continuation and interior point methods rely on a monolithic solution of the discrete energy minimization problem using the trust region minimization algorithm. Whereas the continuation method is based on successive smooth approximations to the minimization problem, the interior point method uses techniques developed for convex second-order cone
programming. The number of load steps was different in each algorithm and was chosen by trial and error to optimize convergece of the iterations. Significantly smaller steps were required in the block-CD algorithm to obtain convergent iterations. 
The ADMM outperforms all previous algorithms in this quasi-static problem. This is due not only to savings because of relatively large load increments, but also the need for fewer iterations for convergence at each step. 

\begin{figure}
\begin{center}
\includegraphics[max width=0.65\linewidth]{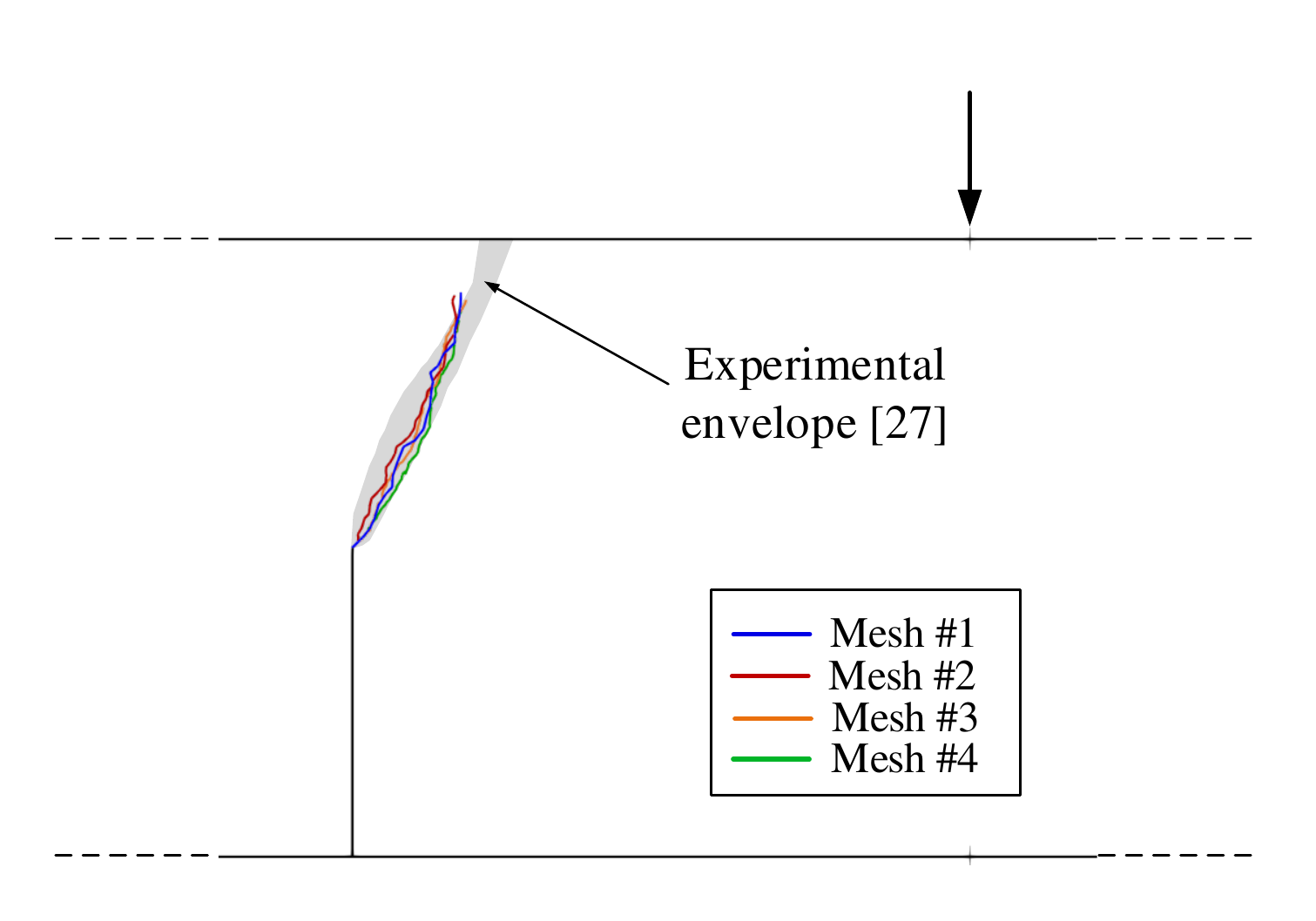}
\caption{Crack paths obtained using different meshes for the single-edge notched beam.}
\label{fig:galvez-paths-meshes}
\end{center}
\end{figure}

Figure \ref{fig:galvez-load-disp-meshes} also shows force-displacement curves obtained using four different meshes. The meshes consist of 5684, 14663, 44773, and 65451 elements, respectively. Furthermore, Figure \ref{fig:galvez-paths-meshes} shows the crack paths.

\begin{table}
\centering
\caption{Comparison of run-times}

{\small
\begin{tabular}{ c c c c c}
 \hline
 \rowcolor{Gray} Algorithm  & Load steps & CPU & Run-time \\ 
 \hline
 Block-CD & 4000 & 3.3 GHz Core i7 5820K & 9.6 hr \\      
 \hline
 Continuation &  500 & 3.3 GHz Core i7 5820K & 13 hr \\  
 \hline
 Interior point \cite{vavasis2020second} &  26 & 2.6 GHz Xeon E5-2690 &  17 hr \\
 \hline
 ADMM (without extap.) & 200 & 3.3 GHz Core i7 5820K &  21.7 min \\
 \hline
 ADMM (with extrap.) &  200 & 3.3 GHz Core i7 5820K & 3.1 min \\
 \hline
\end{tabular}}
\label{tab:comparison}
\end{table}


\begin{figure}

\begin{subfigure}{0.6\textwidth}
\begin{center}
{\includegraphics[max width=0.9\linewidth]{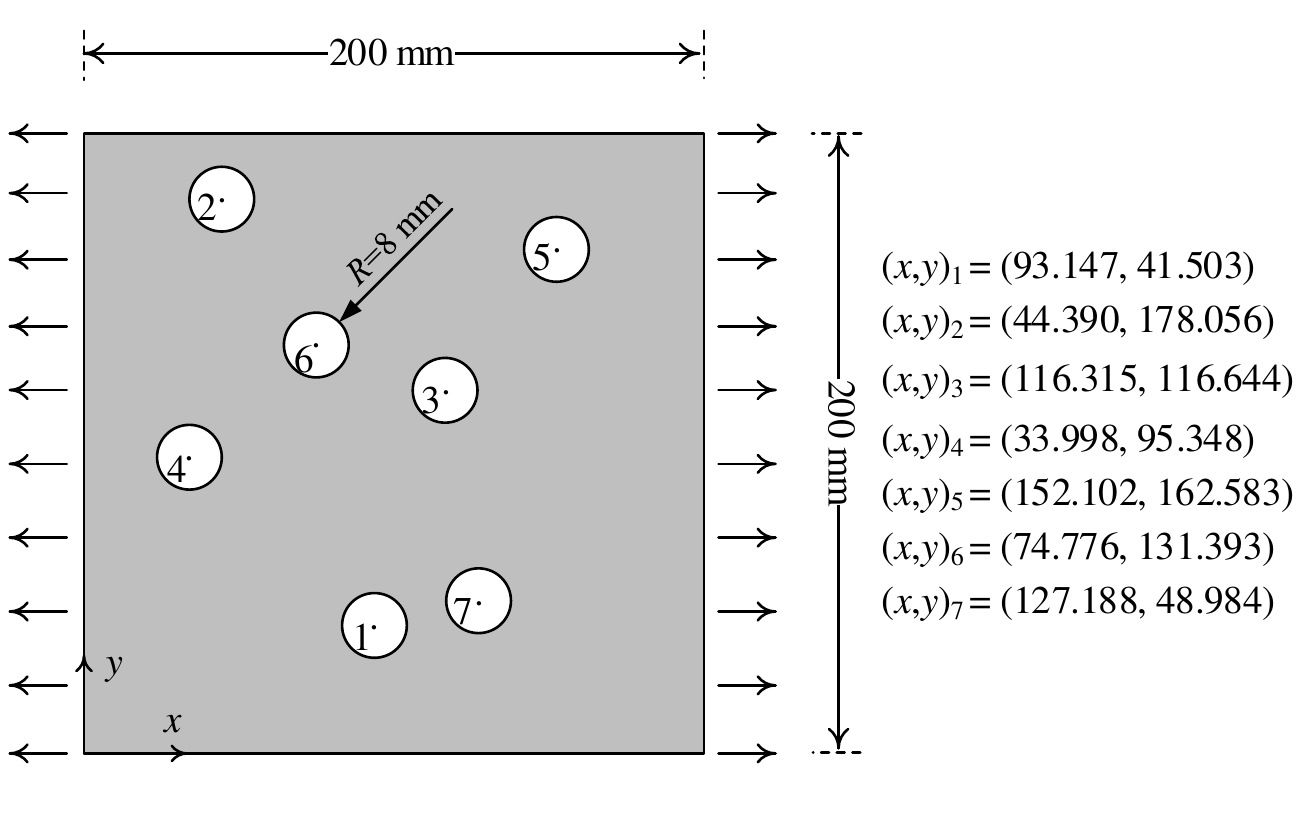}}
\caption{}
\end{center}
\end{subfigure}
\hfill
\begin{subfigure}{0.33\textwidth}
\begin{center}
{\includegraphics[max width=0.9\linewidth]{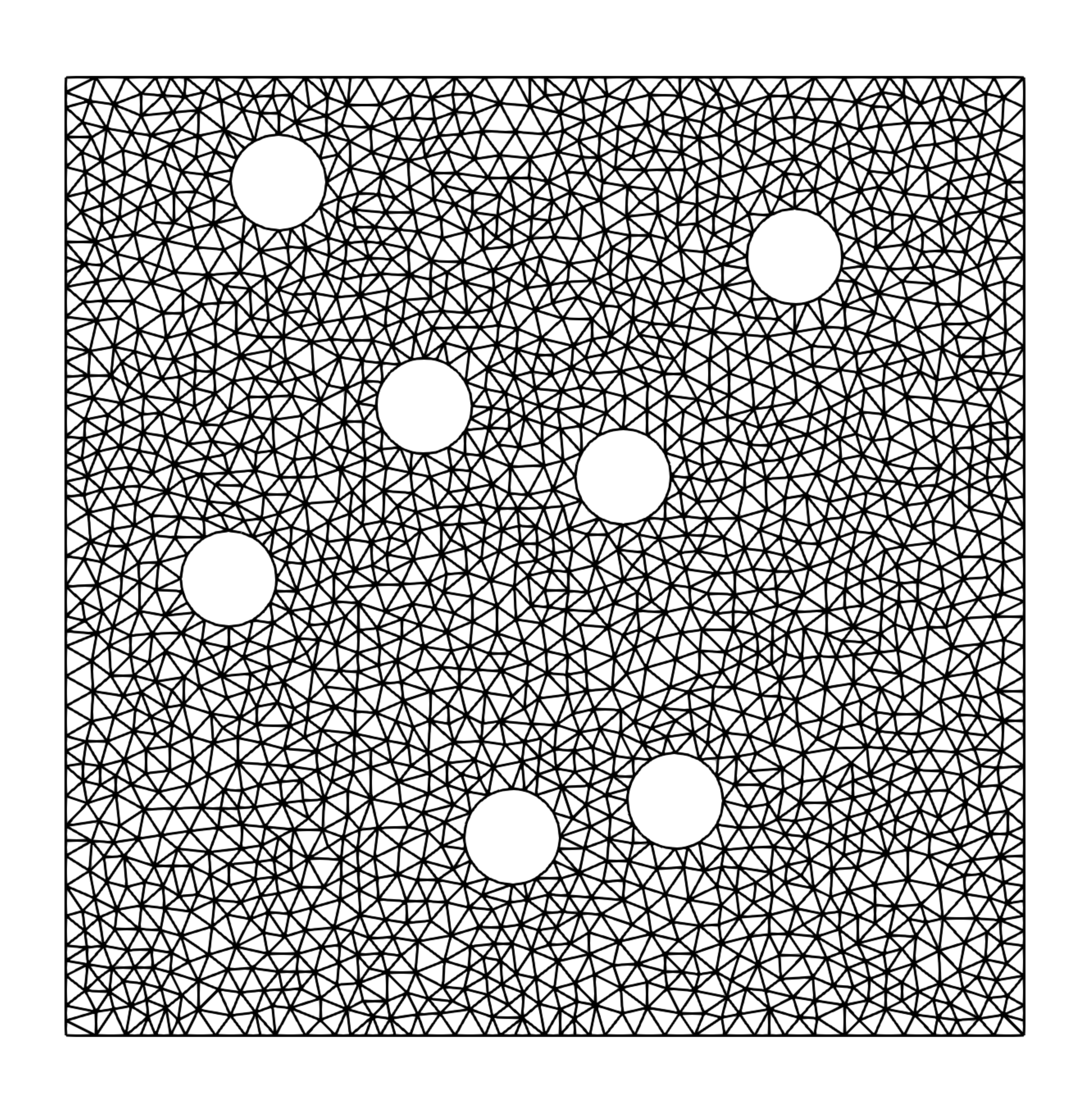}}
\caption{}
\end{center}
\end{subfigure}

\caption{Test problem; a) geometry and boundary conditions and, b) typical mesh used in the simulations containing 3913 elements.}
\label{fig:test_problem}
\end{figure}

\section{Test problem}
\label{test}
Our purpose in this section is to further test the algorithm for performance and sensitivity to algorithmic parameters. 
To this end, 
a test problem was devised, 
which provides insight into the model's performance under different conditions and establishes the effect of different parameter choices. The geometry is based on a material with microstructure in the form of randomly embedded pores. The geometry can be seen in Figure \ref{fig:test_problem},
in which the pores are randomly generated but required to be a certain distance apart. In these tests, a lateral uniform displacement is incrementally applied on a vertical edge and the sum of lateral forces on the edge nodes is measured at each step. The displacement and force are converted to ``stress'' and ``strain'' by dividing by the width and height of the geometry respectively. A pinwheel isoperimetric mesh \cite{pinwheel} was used for each of the tests to reduce the mesh dependence of the solutions. The default mesh has both width and length set to 200 and consists of 7788 elements, $\sigma_c$ is set to 3.0 and $\delta_c$ to $0.02287$. The default values for $c_{dual}$ and $c_{primal}$ are set to 0.01 and about 1/300 of the maximum local stress during fracture (based on $\sigma_c$). Extrapolation and irreversibility are enabled by default, and the displacement is applied in 200 uniform steps between 0 and  0.029731 (slightly larger than $\delta_c$).

A parameter $\alpha$ is introduced to ensure that the convexity condition given by equation \eqref{eqn:rhoCriteria} is satisfied when other parameters, including the mesh resolution, are changed. $\rho$ is computed from $\rho = \alpha \times mean(a_i)\times\frac{\sigma_c}{\delta_c}$. The mean is chosen instead of the maximum in order to make $\rho$ less sensitive to small irregularities in the mesh but the exact criterion of \eqref{eqn:rhoCriteria} is also verified at each Gauss point. The optimization process is dependent on the magnitude of $\rho$ relative to other parameters so, by setting $\alpha$ constant, some of this dependence is reduced. The default value for $\alpha$ used in these tests is 100, however a range of values is also explored.

\subsection{Convergence of the optimization problem}

\begin{figure}[t!]
\centering
\begin{minipage}{0.49\textwidth}
    \centering
    \resizebox{\linewidth}{!}{
    \begin{tikzpicture}
    \begin{axis}[
        xlabel={Strain $\%$},
        ylabel={Stress [MPa]},
        xmin=0, xmax=0.0175,
        ymin=0, ymax=2.5,
        xtick={0,0.005,0.01,0.015},
        ytick={0,0.5,1.0,1.5,2.0,2.5},
        legend pos=north east,
        /tikz/font=\small,
        legend style={font=\small}
    ]
    \addplot[color=blue, mark=o, mark repeat=5]
        table {./Data/extraplStress_true.dat};
        \addlegendentry{With extrapolation};

    \addplot[color=red, mark=asterisk, mark repeat=5]
        table {./Data/extraplStress_false.dat};
        \addlegendentry{Without extrapolation};
    
    \end{axis}
    \end{tikzpicture}
    }
    \subcaption{}
    \label{fig:extrap-a}
\end{minipage}    
\hfill
\begin{minipage}{0.49\textwidth}
    \centering
    \resizebox{\linewidth}{!}{
    \begin{tikzpicture}
    \begin{axis}[
        xlabel={Step},
        ylabel={Number of optimization iterations},
        xmin=0, xmax=210,
        ymin=0, ymax=250,
        xtick={0,50,100,150,200},
        ytick={0,50,100,150,200, 250},
        legend pos=north east,
        /tikz/font=\small,
        legend style={font=\small}
    ]
    \addplot[color=blue, mark=o, mark repeat=2]
        table {./Data/extrapIters_true.dat};
        \addlegendentry{With extrapolation};

    \addplot[color=red, mark=asterisk, mark repeat=2]
        table {./Data/extrapIters_false.dat};
        \addlegendentry{Without extrapolation};
    
    \end{axis}
    \end{tikzpicture}
    }
    \subcaption{}
    \label{fig:extrap-b}
\end{minipage}    
\caption{a) Stress-strain and b) number of optimization iterations per load step with extrapolation enabled and disabled. Strain increases linearly with load step. The stress-strain curves computed using the two methods are visually indistinguishable. 
For steps, during which a significant number of new cracks are formed (around step 80), the update is not easily predicted using linear extrapolation and a similar number of optimization iterations are needed for the two methods. For most other steps, 
with a linear prediction,
about 90\% of the computation is skipped.}
\label{fig:extrap}
\end{figure}
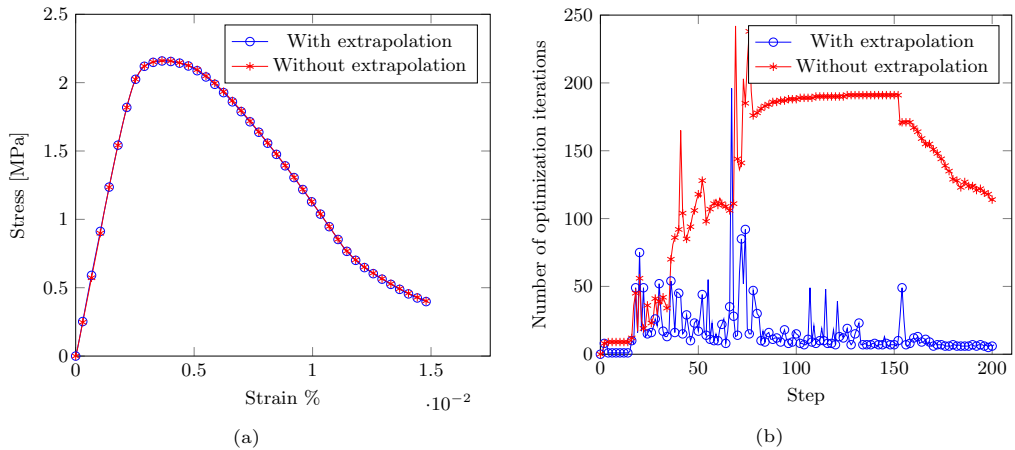

The number of ADMM iterations required for convergence at each load step (as well as the computation time for each iteration) is recorded using the geometry shown in Figure \ref{fig:test_problem}. The average stress - average strain curve and the number of optimization iterations for each load step are shown with and without the extrapolation method in Figure \ref{fig:extrap}. The stress-strain curve is visually indistinguishable in the two cases. The number of iterations is reduced from 26975 to 3327 when extrapolation is activated.

\begin{figure}

\centering
\begin{subfigure}{.3\textwidth}
  \centering
  \includegraphics[width=\linewidth]{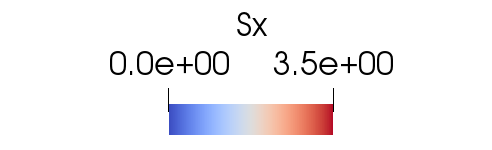}
\end{subfigure}
\vfill
\begin{subfigure}{.28\textwidth}
  \centering
  \includegraphics[width=0.9\linewidth]{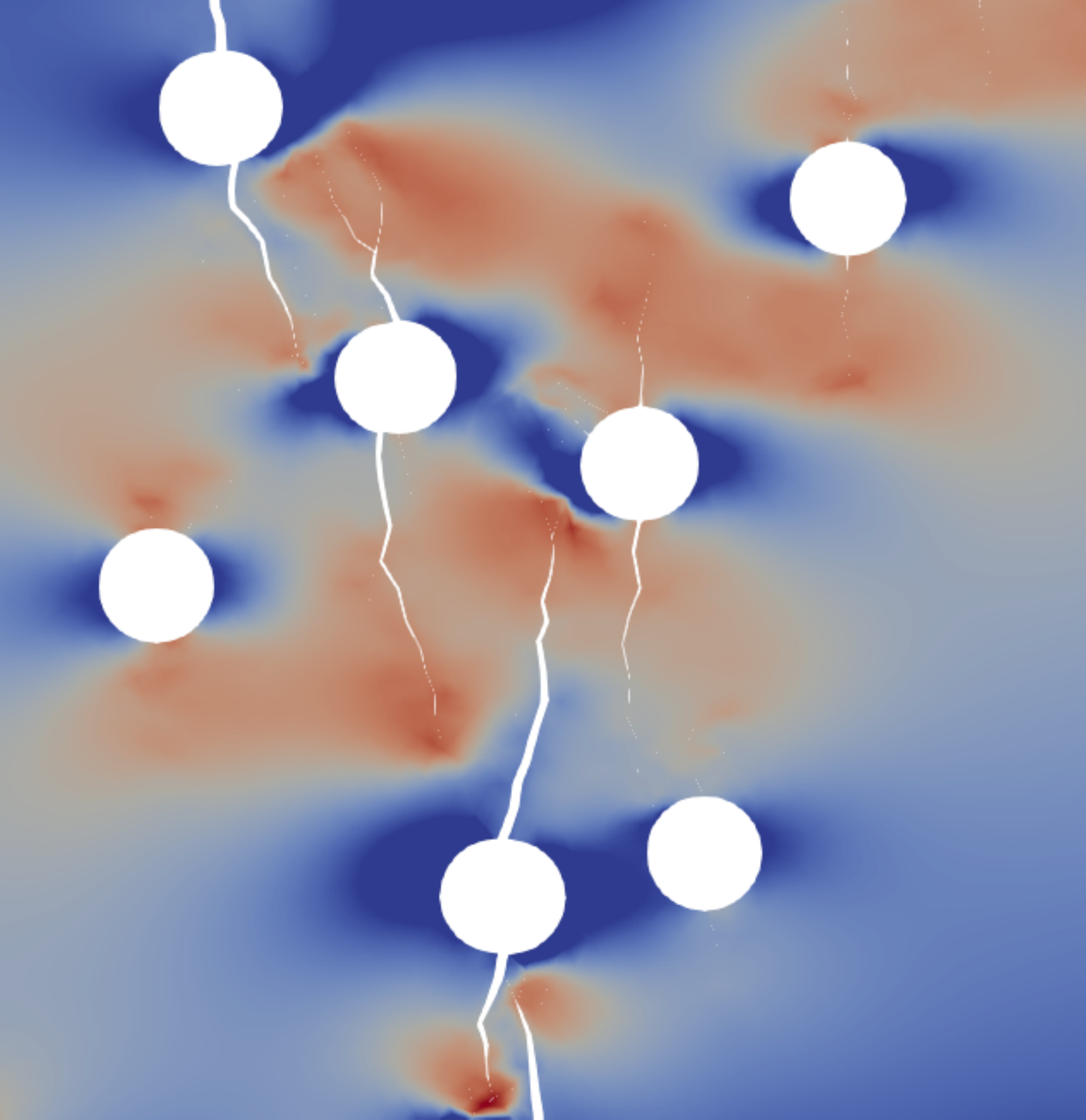}
  \subcaption{3913 elements}
\end{subfigure}
\begin{subfigure}{.28\textwidth}
  \centering
  \includegraphics[width=0.9\linewidth]{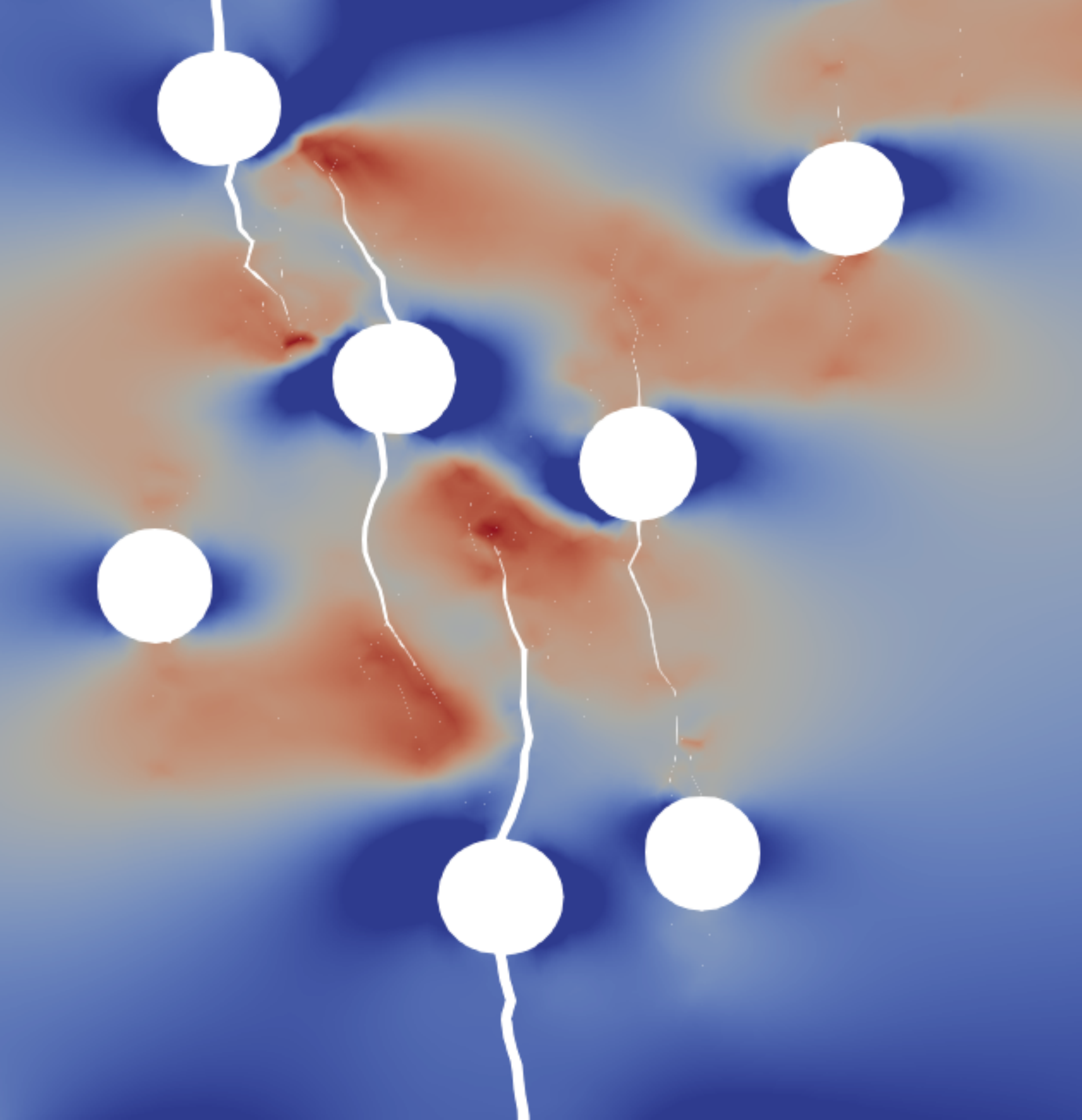}
  \subcaption{5430 elements}
\end{subfigure}

\begin{subfigure}{.28\textwidth}
  \centering
  \includegraphics[width=0.9\linewidth]{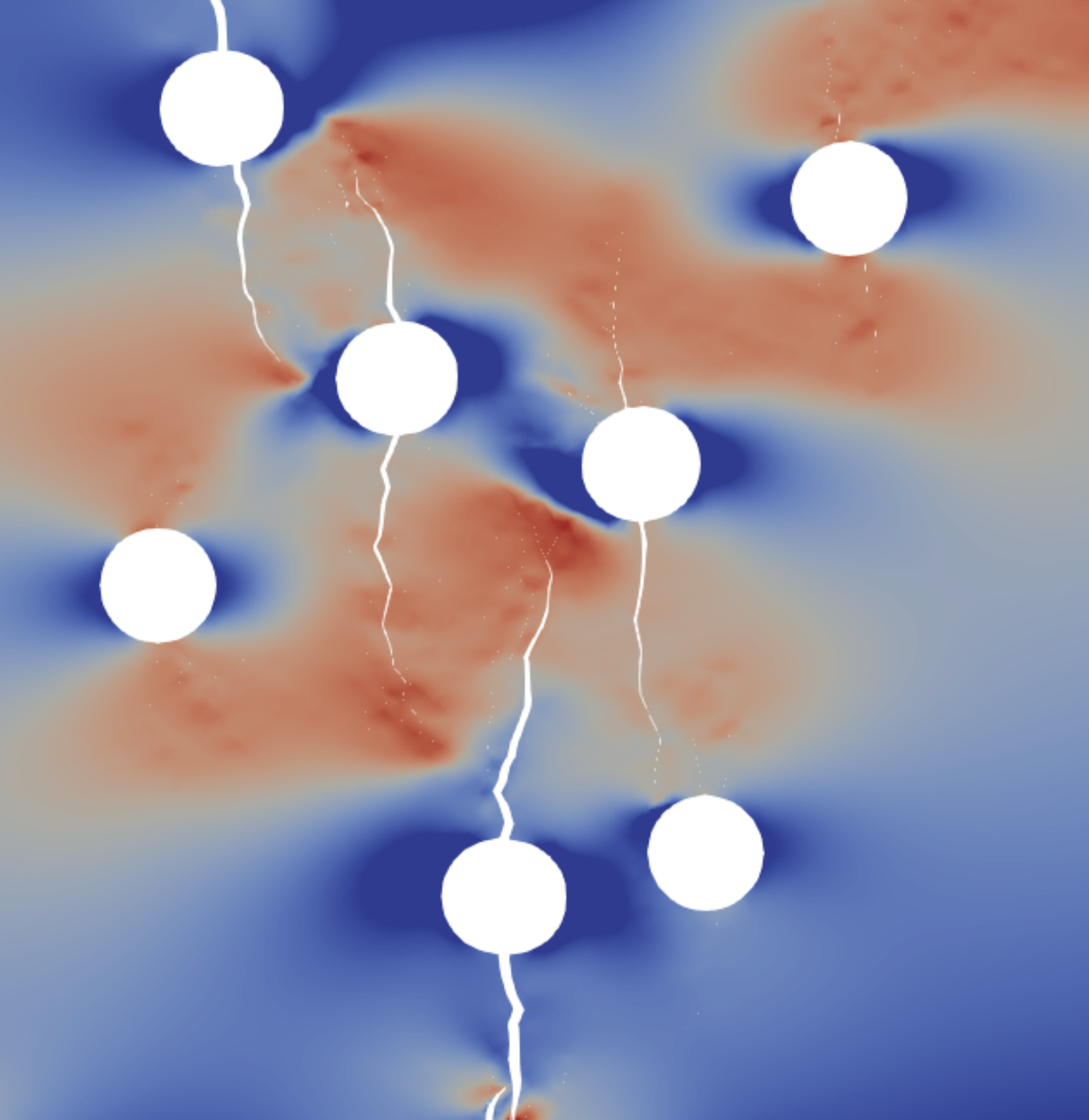}
  \subcaption{11910 elements}
\end{subfigure}
\vfill
\vspace{3pt}
\begin{subfigure}{.28\textwidth}
  \centering
  \includegraphics[width=0.9\linewidth]{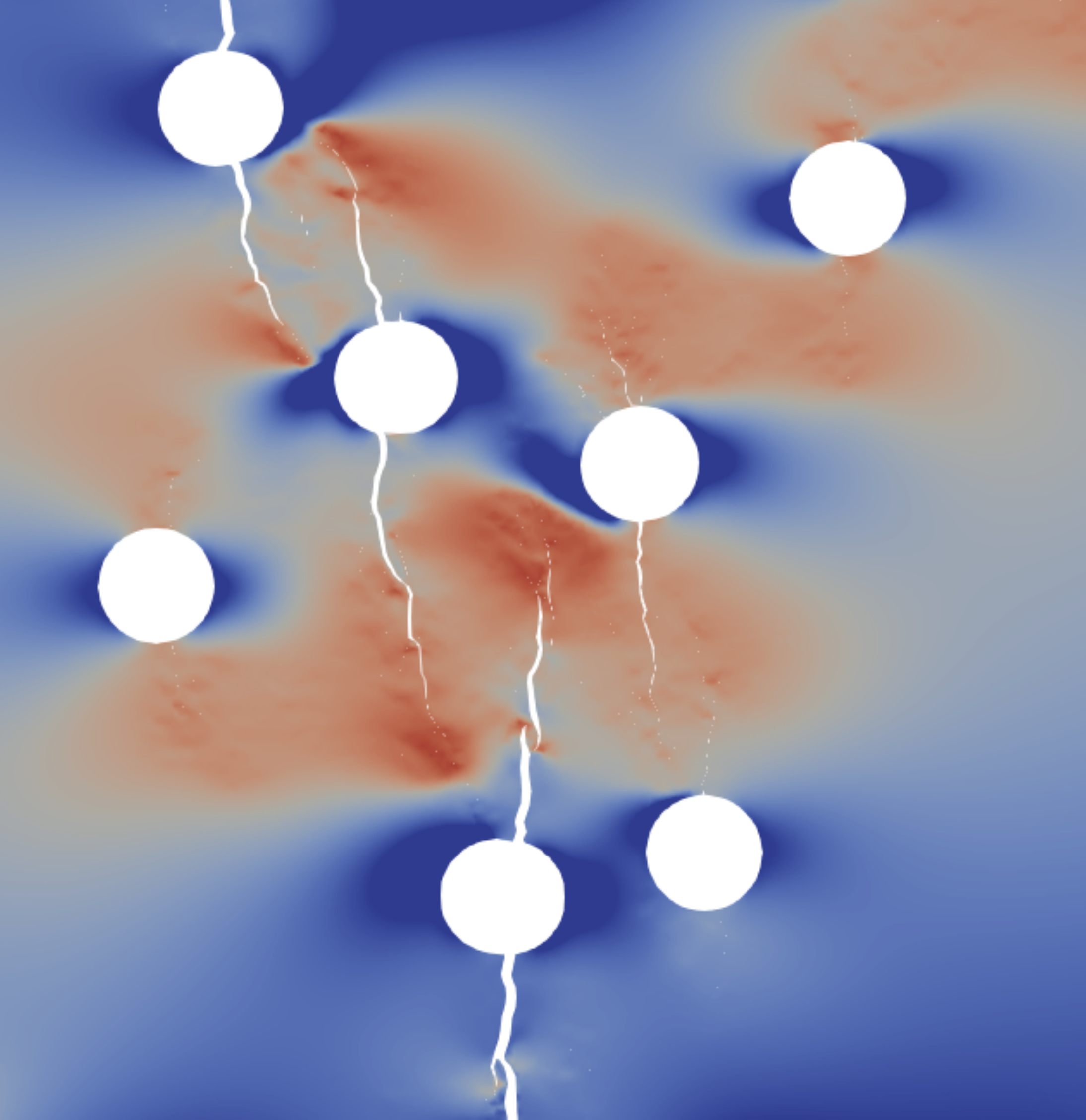}
  \subcaption{21621 elements}
\end{subfigure}
\begin{subfigure}{.28\textwidth}
  \centering
  \includegraphics[width=0.9\linewidth]{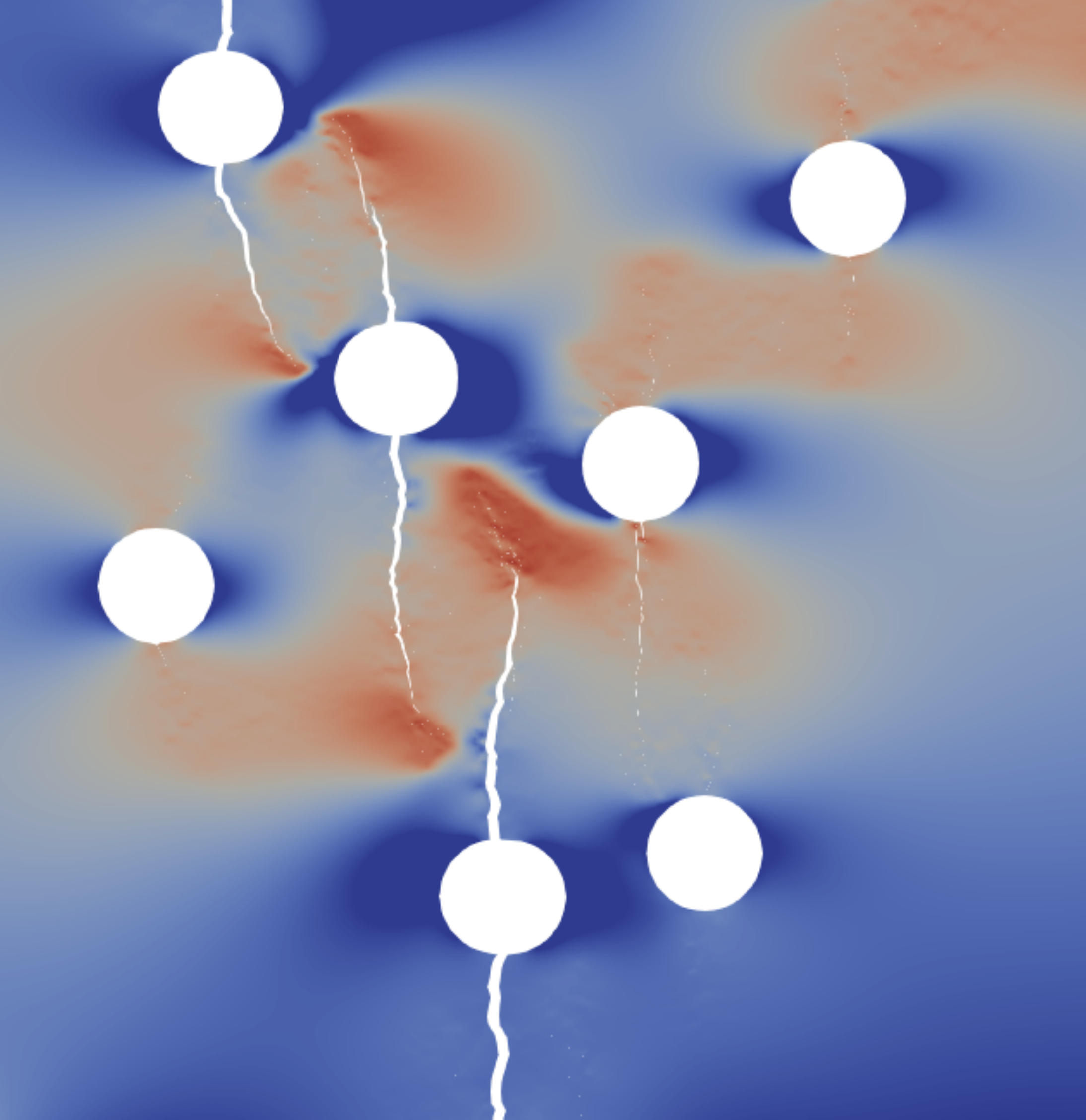}
  \subcaption{47930 elements}
\end{subfigure}

\caption{Crack paths obtained using different meshes for the test problem. Although the paths are consistent, differences in the state of stress and strain in the neighboring elements lead to variations in the force-displacement curve, see Fig. \ref{fig:parameters-c}.}
\label{fig:fracPath}
\end{figure}

Initially, the problem is completely elastic and the number of optimization iterations is small because none of the interfaces open. In this regime, the problem is equivalent to linear elastic FEA, with the added complexity that node displacements from neighbouring elements are constrained to be equal instead of shared. As the strain increases, stress concentrates in certain regions (typically near the edges of pores) and 
reaches $\sigma_c$. Some interfaces open slightly, which increases the number of optimization iterations around step 50. Around step 75, the cracks seen in Figure \ref{fig:fracPath} begin to open substantially and the stress reaches a plateau. This results in the largest change in configuration and a higher number of optimization iterations (both with and without extrapolation). After this peak, the number of iterations required when extrapolation is disabled are similar for neighbouring steps. This is the original motivation for the extrapolation method: if the number of iterations is nearly identical, the problems being solved at each step are likely similar. Therefore there should be useful information available from the optimization process of previous steps, and a linear extrapolation is the simplest kind to try first. With extrapolation enabled, the number of iterations is consistently quite low after the peak, supporting the hypothesis that there is useful information in the trend of previous solutions. \par

Another interesting feature of Figure \ref{fig:extrap} is the drop in the number of iterations without  extrapolation around step 150 at the same time that the number of iterations with extrapolation increases slightly. This is the point when some cracks open completely ($\delta > \delta_c$), and the slope of the stress-strain curve changes (around strain 0.012\%). The iterations without extrapolation decrease because the complexity of the problem is reduced when some cracks open completely. However, the iterations with extrapolation momentarily increase because the initial guess is less accurate due to the deviation from the previous pattern.

\subsection{Parameter sensitivity} 
\label{ParamSensitivity}

The main parameters of the optimization process are $\rho$,
$c_{\rm primal}$, and $c_{\rm dual}$. These impact the priority given to satisfying the equality constraint during intermediate optimization iterations as well as the stopping criteria. 
Another setting
decides whether extrapolation is enabled. 
These change the initial guess at each optimization step and the way the potential energy is modified as the solution progresses. The problem definition is dependent on the mesh resolution and on the load step resolution. It is useful to analyze both the sensitivity of solutions, and how changes in these parameters affect runtime.

\begin{figure}[t!]
\centering
\begin{minipage}{0.49\textwidth}
    \centering
    \footnotesize
    \resizebox{\linewidth}{!}{
    \begin{tikzpicture}
    \begin{axis}[
        xlabel={Strain $\%$},
        ylabel={Stress [MPa]},
        xmin=0, xmax=0.0175,
        ymin=0, ymax=2.5,
        xtick={0,0.005,0.01,0.015},
        ytick={0,0.5,1.0,1.5,2.0,2.5},
        legend pos=north east,
        /tikz/font=\small,
        legend style={font=\small}
    ]
    \addplot[color=blue, mark=o, mark repeat=5]
        table {./Data/stress_alpha20.dat};
        \addlegendentry{$\alpha$=20};

    \addplot[color=red, mark=asterisk, mark repeat=5]
        table {./Data/stress_alpha60.dat};
        \addlegendentry{$\alpha$=60};

    \addplot[color=orange, mark=x, mark repeat=5]
        table {./Data/stress_alpha100.dat};
        \addlegendentry{$\alpha$=100};
        
    \addplot[color=teal, mark=triangle, mark repeat=5]
        table {./Data/stress_alpha180.dat};
        \addlegendentry{$\alpha$=180};        
    
    \end{axis}
    \end{tikzpicture}
    }
    \subcaption{}
    \label{fig:parameters-a}
\end{minipage}    
\begin{minipage}{0.49\textwidth}
    \centering
    \footnotesize
    \resizebox{\linewidth}{!}{
    \begin{tikzpicture}
    \begin{axis}[
        xlabel={Strain $\%$},
        ylabel={Stress [MPa]},
        xmin=0, xmax=0.0175,
        ymin=0, ymax=2.5,
        xtick={0,0.005,0.01,0.015},
        ytick={0,0.5,1.0,1.5,2.0,2.5},
        legend pos=north east,
        /tikz/font=\small,
        legend style={font=\small}
    ]
    \addplot[color=blue, mark=o, mark repeat=1]
        table {./Data/stress_load25.dat};
        \addlegendentry{27 steps};

    \addplot[color=red, mark=asterisk, mark repeat=2]
        table {./Data/stress_load100.dat};
        \addlegendentry{102 steps};

    \addplot[color=orange, mark=x, mark repeat=8]
        table {./Data/stress_load400.dat};
        \addlegendentry{402 steps};
        
    \addplot[color=teal, mark=triangle, mark repeat=32]
        table {./Data/stress_load1600.dat};
        \addlegendentry{1602 steps};     
    
    \end{axis}
    \end{tikzpicture}
    }
    \subcaption{}
    \label{fig:parameters-b}
\end{minipage}    
\vfill
\begin{minipage}{0.49\textwidth}
    \centering
    \footnotesize
    \resizebox{\linewidth}{!}{
    \begin{tikzpicture}
    \begin{axis}[
        xlabel={Strain $\%$},
        ylabel={Stress [MPa]},
        xmin=0, xmax=0.0175,
        ymin=0, ymax=2.5,
        xtick={0,0.005,0.01,0.015},
        ytick={0,0.5,1.0,1.5,2.0,2.5},
        legend pos=north east,
        /tikz/font=\small,
        legend style={font=\footnotesize}
    ]
    \addplot[color=blue, mark=o, mark repeat=5]
        table {./Data/stress_elem3913.dat};
        \addlegendentry{3913 elements};

    \addplot[color=red, mark=asterisk, mark repeat=5]
        table {./Data/stress_elem5430.dat};
        \addlegendentry{5430 elements};

    \addplot[color=orange, mark=x, mark repeat=5]
        table {./Data/stress_elem11910.dat};
        \addlegendentry{11910 elements};
        
    \addplot[color=magenta, mark=triangle, mark repeat=5]
        table {./Data/stress_elem21624.dat};
        \addlegendentry{21624 elements};     
        
    \addplot[color=teal, mark=diamond, mark repeat=5]
        table {./Data/stress_elem47930.dat};
        \addlegendentry{47930 elements};            
    
    \end{axis}
    \end{tikzpicture}
    }
    \subcaption{}
    \label{fig:parameters-c}
\end{minipage}    
\begin{minipage}{0.49\textwidth}
    \centering
    \footnotesize
    \resizebox{\linewidth}{!}{
    \begin{tikzpicture}
    \begin{axis}[
        xlabel={Strain $\%$},
        ylabel={Stress [MPa]},
        xmin=0, xmax=0.0175,
        ymin=0, ymax=2.5,
        xtick={0,0.005,0.01,0.015},
        ytick={0,0.5,1.0,1.5,2.0,2.5},
        legend pos=north east,
        /tikz/font=\small,
        legend style={font=\footnotesize}
    ]
    \addplot[color=blue, mark=o, mark repeat=5]
        table {./Data/stress_tol1024.dat};
        \addlegendentry{PressTol = 1.024};

    \addplot[color=red, mark=asterisk, mark repeat=5]
        table {./Data/stress_tol256.dat};
        \addlegendentry{PressTol = 0.256};

    \addplot[color=orange, mark=x, mark repeat=5]
        table {./Data/stress_tol016.dat};
        \addlegendentry{PressTol = 0.016};
        
    \addplot[color=magenta, mark=triangle, mark repeat=5]
        table {./Data/stress_tol001.dat};
        \addlegendentry{PressTol = 0.001};     
        
    \addplot[color=teal, mark=diamond, mark repeat=5]
        table {./Data/stress_tol00025.dat};
        \addlegendentry{PressTol = 0.00025};            
    
    \end{axis}
    \end{tikzpicture}
    }
    \subcaption{}
    \label{fig:parameters-d}
\end{minipage}    

\caption{Stress-strain curves as lateral strain is applied for various parameter ranges. The parameters varied are a) $\alpha$, b) number of load steps, c) number of mesh elements, and d) convergence threshold. In (a) and (b), all solutions are nearly identical. For the mesh size, 
The apparent lack of convergence with mesh size is believed to be due to the existence of closely spaced minima.
As the stopping criteria decrease the stress-strain curves converge.}
\label{fig:parameters}
\end{figure}
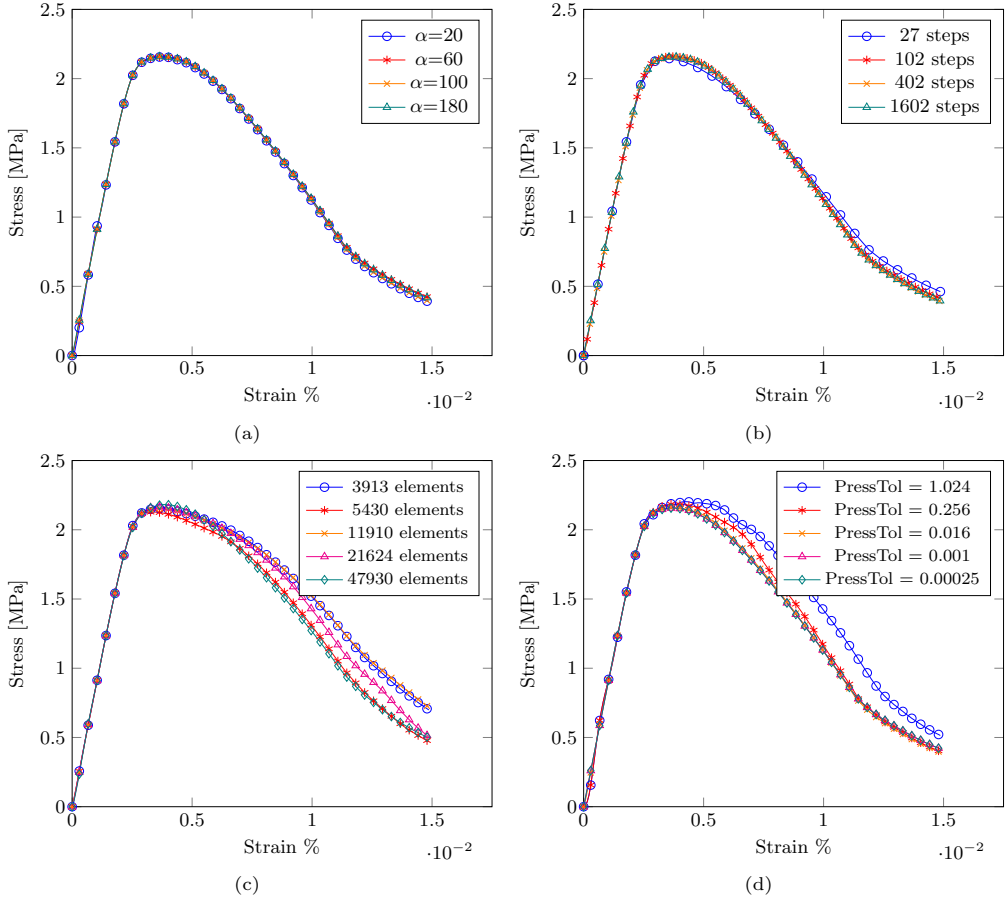

Figure \ref{fig:parameters} shows stress-strain curves for parameter sweeps of $\alpha$, $c_{primal} \text{ \& } c_{dual}$, mesh resolution, and load resolution. The number of optimization iterations needed to find the solution for each choice of parameters is discussed in detail in Section \ref{time}. 

The results are almost completely insensitive to $\alpha$ (while kept a safe distance away from the convexity condition near $\alpha=1$); see Figure \ref{fig:parameters-a}. Near the right side of the plot, there is a slight deviation between high and low $\alpha$ values. This is likely due to slightly different local minima being found based on how much the optimization algorithm relaxes the equality constraints during intermediate steps. 

The effect of displacement increment size is considered in Figure \ref{fig:parameters-b}. In general, larger steps may allow the state to escape a local minimum. For this example, the optimization finds the same local minimum for all displacement increment sizes except when the number of steps is set to the lowest tested value, 27. With 27 steps the progression is less smooth and deviates slightly from the trajectory obtained with other step sizes.  The convergence with decreasing load step size is expected,  being the quasistatic equivalent of convergence with time step, which is known to hold for time-continuous models but not for extrinsic cohesive models which do not possess the time-continuous property \cite{timecont, contDG}.

We test a moderate range of mesh resolutions to estimate mesh sensitivity rather than mesh convergence. 
Figure \ref{fig:parameters-c} shows mesh dependence in the falling branch of the force-displacement curve in the sense that solutions do not seem to become closer with mesh refinement. The crack paths obtained using these meshes are shown in Fig. \ref{fig:fracPath}. In particular, note the closeness of the 3913 and 11910 elements solutions and that between solutions obtained with 5430 and 47930 elements. This clustering, and the fact that no observable mesh dependence typically occurs with meshes possessing the properties of isotropy and isoperimetry \cite{isotropicmesh}, e.g., pinwheel meshes \cite{pinwheel},  but also with general unstructured meshes, leads us to believe that closely spaced minima exist in problems with this type of microstructure such that a variation in the mesh will pick up different solutions corresponding to these minima. This problem-specific variability - some would call it a ``data'' (as opposed to algorithmic) instability - is probably random in nature and will certainly be picked up by any algorithm. Since this type of microstructure is typical of many multiscale analyses, a stochastic analysis of representative volumes seems necessary.  A similar conclusion was reached by Al-Ostaz and Jasiuk \cite{al1997crack} in the study of  perforated epoxy and aluminum sheets. They performed Ansys$\mbox{}^{\mathrm{TM}}$ and spring-network fracture simulations, in which they obtained mesh dependence of the crack paths and, more importantly, they experimentally observed variability of the force-displacement curve. 

Finally, the effect of the tolerance used in the stopping criteria is shown in Figure \ref{fig:parameters-d}. Both primal and dual pressures were required to be less than the `PressTol' parameter for convergence. As the criteria become stricter (lower 
tolerance), the solutions converge to a single curve. For tolerance values of 0.016 and below, the results are indistinguishable.

\subsection{Time complexity}
\label{time}

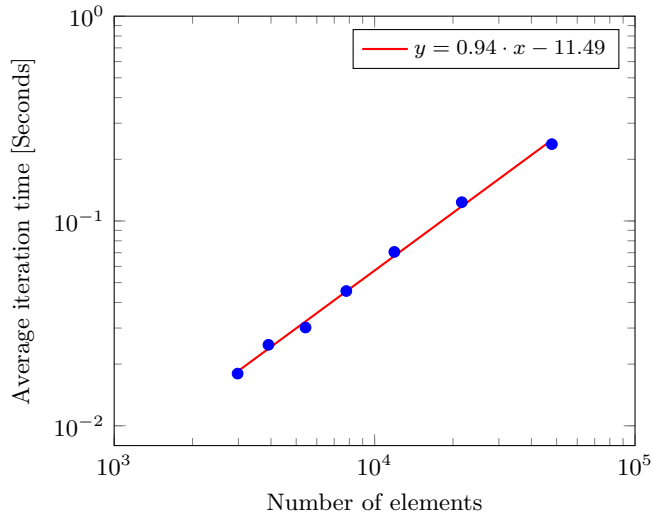
\begin{figure}[ht]
\centering
    \begin{tikzpicture} 
    \begin{axis}[
        xmode=log,
        ymode=log,
        xlabel={Number of elements},
        ylabel={Average iteration time [Seconds]},
        xmin=1000, xmax=100000,
        ymin=0.008, ymax=1.0,
        xtick={ 1000, 10000, 100000},
        ytick={0.01,0.1,1.0},
        legend pos=north east,
        /tikz/font=\small,
        legend style={font=\footnotesize}
    ]

    \addplot [thick, red] table[
        y={create col/linear regression={y=Y}}
    ] 
    {./Data/IterTime.dat};    
    \addlegendentry{
        $ y =
        \pgfmathprintnumber{\pgfplotstableregressiona}
        \cdot x
        \pgfmathprintnumber[print sign]{\pgfplotstableregressionb}$
    };     
    \addplot[color=blue, mark=*, only marks]
        table {./Data/IterTime.dat};
   


        



    \end{axis}
    \end{tikzpicture}
\caption{Average number of seconds per optimization iteration for various numbers of mesh elements. The relationship is nearly linear, as expected.}
\label{fig:iterTime}
\end{figure}

The run time of the model is analyzed for the previous set of experiments. 
We consider the average time it takes to complete each optimization iteration and the total number of iterations separately.  The tests were performed using a single core on an Intel Xeon CPU. Figure \ref{fig:iterTime} shows the average time needed for each iteration as a function of the problem size. Since the bulk of the code is composed of sparse matrix-vector multiplications and single-pass loops, the iteration time is expected to be nearly linear with problem size. The fitted curve in the figure confirms this, with an exponent close to 0.94. The per-iteration computation time is independent of  other parameters of interest.

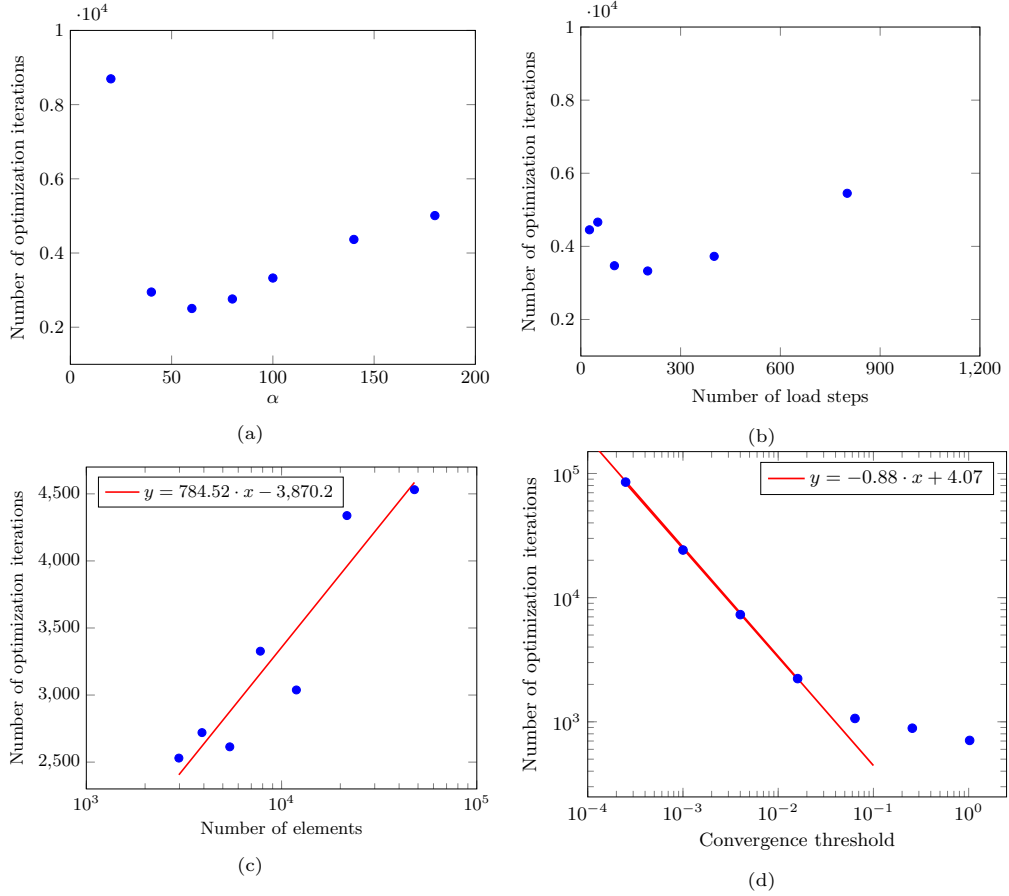
\begin{figure}[t!]
\centering
\begin{minipage}{0.49\textwidth}
    \centering
    \footnotesize
    \resizebox{\linewidth}{!}{
    \begin{tikzpicture}
    \begin{axis}[
        xlabel={$\alpha$},
        ylabel={Number of optimization iterations},
        xmin=0, xmax=200,
        ymin=1000, ymax=10000,
        xtick={0,50,100,150,200},
        ytick={2000,4000,6000,8000,10000},
        legend pos=north east,
        /tikz/font=\small,
        legend style={font=\small}
    ]
    \addplot[color=blue, mark=*, only marks]
        table {./Data/iters_Rho.dat};
    
    \end{axis}
    \end{tikzpicture}
    }
    \subcaption{}
    \label{fig:time-parameters-a}
\end{minipage}    
\begin{minipage}{0.49\textwidth}
    \centering
    \footnotesize
    \resizebox{\linewidth}{!}{
    \begin{tikzpicture}
    \begin{axis}[
        xlabel={Number of load steps},
        ylabel={Number of optimization iterations},
        xmin=0, xmax=1200,
        ymin=1000, ymax=10000,
        xtick={0,300,600,900,1200},
        ytick={2000,4000,6000,8000,10000},
        legend pos=north east,
        /tikz/font=\small,
        legend style={font=\small}
    ]
    \addplot[color=blue, mark=*, only marks]
        table {./Data/iters_loads.dat};
    
    \end{axis}
    \end{tikzpicture}
    }
    \subcaption{}
    \label{fig:time-parameters-b}
\end{minipage}    
\vfill
\begin{minipage}{0.49\textwidth}
    \centering
    \footnotesize
    \resizebox{\linewidth}{!}{
    \begin{tikzpicture}
    \begin{axis}[
        xmode=log,
        xlabel={Number of elements},
        ylabel={Number of optimization iterations},
        xmin=1000, xmax=100000,
        ymin=2300, ymax=4700,
        xtick={100, 1000,10000,100000},
        ytick={2500,3000,3500,4000,4500,5000},
        legend pos=north west,
        /tikz/font=\small,
        legend style={font=\small}
    ]
    \addplot [thick, red] table[
        y={create col/linear regression={y=Y}}
    ] 
    {./Data/iters_mesh.dat};       
    \addlegendentry{
        $ y =
        \pgfmathprintnumber{\pgfplotstableregressiona}
        \cdot x
        \pgfmathprintnumber[print sign]{\pgfplotstableregressionb}$
    };       
    \addplot[color=blue, mark=*, only marks]
        table {./Data/iters_mesh.dat};

    \end{axis}
    \end{tikzpicture}
    }
    \subcaption{}
    \label{fig:time-parameters-c}
\end{minipage}    
\begin{minipage}{0.49\textwidth}
    \centering
    \footnotesize
    \resizebox{\linewidth}{!}{
    \begin{tikzpicture}
    \begin{axis}[
        xmode=log,
        ymode=log,
        xlabel={Convergence threshold},
        ylabel={Number of optimization iterations},
        xmin=0.0001, xmax=2.5,
        ymin=0, ymax=150000,
        xtick={0.0001, 0.001, 0.01, 0.1, 1.0},
        ytick={10, 100,1000,10000,100000},
        legend pos=north east,
        /tikz/font=\small,
        legend style={font=\small}
    ]
    \addplot [thick, red] table[
        y={create col/linear regression={y=Y}}
    ] 
    {./Data/iters_tols_initial.dat};    
        \addlegendentry{
        $ y =
        \pgfmathprintnumber{\pgfplotstableregressiona}
        \cdot x
        \pgfmathprintnumber[print sign]{\pgfplotstableregressionb}$
    };    
    \addplot [domain=0.0001:0.1, thick, red] {exp(-0.88 * ln(x) + 4.07)}; 
    \addplot[color=blue, mark=*, only marks]
        table {./Data/iters_tols.dat};

    \end{axis}
    \end{tikzpicture}
    }
    \subcaption{}
    \label{fig:time-parameters-d}
\end{minipage}    

\caption{Total number of optimization iterations (summed across all load steps) per simulation run for various parameter values. The parameters varied are: top left - $\alpha$, top right - number of load steps, bottom left - number of mesh elements, bottom right - convergence threshold (tolerance); $\alpha$ and the number of load steps both demonstrate a minimum value, while the the number of iterations increases sublinearly with number of elements, and converges almost linearly with the convergence threshold. The  number of iterations as a function of Convergence Threshold starts to diverge from the trend line at about 1000 iterations (5x larger than the theoretical minimum of 1 iteration for each of the 200 load steps).}
\label{fig:time-parameters}
\end{figure}

Figure \ref{fig:time-parameters} shows the total number of optimization iterations needed for each choice of parameters. 
The value of $\rho$ has been shown to be problem dependent and to have a large impact on the convergence rate of ADMM \cite{admmcrit}. If the value is too small, the optimization process 
fails to satisfy the constraints and convergence is slow. If the value is too large, the inner optimization problems 
make small progress toward the solution of the outer optimization problem in each outer iteration. For this test case, the number of (outer) iterations is smallest for 
$\alpha=60$. At $\alpha=20$, the number of iterations increases sharply. 
For small values of $\alpha$, the $\delta$ subproblem stops having a unique solution and also small values of $\alpha$ may allow the process to accidentally escape local minima. For this reason, $\alpha = 100$ is used as the default value instead of $\alpha = 60$. The computation time is similar and there is smaller risk of unintentionally escaping minima. 

The overall pattern of number of optimization iterations versus number of load steps, $N_\text{steps}$, seems to be convex with a minimum around $N_\text{steps}=200$. The $N_\text{steps}=27$ point breaks the trend but, as seen in Figure \ref{fig:parameters}, the solution when $N_\text{steps}=27$  diverges to the point that it may be solving a slightly different problem. With a larger number of steps, each optimization is likely closer to the local minimum and the starting point is better predicted by the extrapolation method. However, as the number of steps increases, the number of individual optimizations also increases. It seems as though there is a fixed cost for each step due to the number of iterations it takes to converge near the solution. 
The location of the minimum is likely sensitive to the convergence tolerance.

The number of iterations is also plotted against the number of mesh elements (again as a log-log plot). 
The overall trend suggests a sublinear relation with exponent about equal to 0.2. This is a slower increase than might otherwise be expected, suggesting that the difficulty in solving the optimization problem is fairly independent of mesh resolution.

The impact of $c_{primal}$ and $c_{dual}$ on number of iterations is plotted (in a log log scale). For values below 0.02, the trend is well fit by a line with slope -0.88. This implies a nearly linear convergence rate in this regime. For values larger than 0.02 the trend line doesn't fit as well, but here the average number of iterations per step is around 5, which is quite close to the minimum limit (1). 

The overall scaling exponent is equal to the sum of the exponents for the change in iteration time and number of iterations. Therefore, the run time scales super linearly (~1 + ~0.2 = 1.2) in number of mesh elements. The convergence tolerance does not impact the iteration time, so the previous (almost linear) result holds. 

The amount of time required to generate a reasonable solution to this problem is about 1.5 minutes (with 5000 elements there are ~3000 steps that take 0.03 seconds each). This is sufficiently fast for the original use-case of the ADMM model, which is large scale data generation for machine learning applications. Based on the roughly linear observed time complexity and some preliminary tests, it is believed that the ADMM model can be used to simulate previously intractable problem sizes in a reasonable time-frame.

\section{Conclusion}

A model was developed to simulate cohesive fracture by using ADMM to minimize a non-smooth potential energy function. Representative convergence paths were analyzed and the sensitivity of predictions to changes in model parameters was tested. An extrapolation method for better initial optimization estimates was developed and shown to significantly improve performance with little impact on results. Finally, the run-time of the model was experimentally shown to increase roughly linearly with problem size and the model was shown to simulate moderately sized problems in seconds. This makes the model well suited for multiple runs within machine learning applications.

\bibliography{admm}
\bibliographystyle{model1-num-names}

\end{document}